\documentclass{article}
\usepackage{stmaryrd}
\usepackage{mathrsfs}
\usepackage[centertags]{amsmath}
\usepackage{amsfonts, dsfont}
\usepackage{amssymb}
\usepackage{amsthm}

\input xy
\xyoption{all}

\theoremstyle{plain}
\newtheorem{thm}{Theorem}[section]
\newtheorem{cor}[thm]{Corollary}
\newtheorem{lem}[thm]{Lemma}
\newtheorem{prop}[thm]{Proposition}
\newtheorem*{thm2}{Theorem}

\theoremstyle{definition}

\newtheorem{remark}[thm]{Remark}
\newtheorem*{ack}{Acknowledgments}

\newcommand{\bl}{\begin{lem}}
\newcommand{\el}{\end{lem}}
\newcommand{\bp}{\begin{prop}}
\newcommand{\ep}{\end{prop}}
\newcommand{\bt}{\begin{thm}}
\newcommand{\et}{\end{thm}}
\newcommand{\bc}{\begin{cor}}
\newcommand{\ec}{\end{cor}}
\newcommand{\br}{\begin{remark}}
\newcommand{\er}{\end{remark}}
\newcommand{\ba}{\begin{array}}
\newcommand{\ea}{\end{array}}
\newcommand{\bpf}{\begin{proof}}
\newcommand{\epf}{\end{proof}}

\newcommand{\R}{\mathds{R}}

\newcommand{\Q}{\mathds{Q}}
\newcommand{\Zp}{\mathds{Z}_{p}}
\newcommand{\Qp}{\mathds{Q}_{p}}
\newcommand{\al}{\alpha}
\newcommand{\Ga}{\Gamma}

\newcommand{\La}{\Lambda}
\newcommand{\la}{\lambda}
\newcommand{\Si}{\Sigma}

\DeclareMathOperator{\Sel}{Sel} \DeclareMathOperator{\Gal}{Gal}
\DeclareMathOperator{\Hom}{Hom} \DeclareMathOperator{\rank}{rank}
\DeclareMathOperator{\Ext}{Ext} \DeclareMathOperator{\Ind}{Ind}

\newcommand{\cyc}{\mathrm{cyc}}
\newcommand{\cts}{\mathrm{cts}}

\newcommand{\Op}{\mathcal{O}}

\newcommand{\M}{\mathfrak{M}}
\newcommand{\A}{\mathcal{A}}

\newcommand{\ot}{\otimes}
\newcommand{\ilim}{\displaystyle \mathop{\varinjlim}\limits}
\newcommand{\plim}{\displaystyle \mathop{\varprojlim}\limits}
\newcommand{\im}{\mathrm{im}\,}

\newcommand{\lra}{\longrightarrow}
\newcommand{\tha}{\twoheadrightarrow}

\newcommand{\ps}[1]{\llbracket #1 \rrbracket}

\begin{document}
\title{$\M_H(G)$-property and congruence of Galois representations }
\author{Meng Fai Lim\footnote{School of Mathematics and Statistics $\&$ Hubei Key Laboratory of Mathematical Sciences,
Central China Normal University, Wuhan, 430079, P.R.China.
 E-mail: \texttt{limmf@mail.ccnu.edu.cn}}}
\date{}
\maketitle

\begin{abstract}
{\footnotesize  In this paper, we study the Selmer groups of two congruent
Galois representations over an admissible $p$-adic Lie extension. We
will show that under appropriate congruence condition, if the dual
Selmer group of one satisfies the $\M_H(G)$-property, so will the
other. In the event that the $\M_H(G)$-property holds, and assuming
certain further hypothesis on the decomposition of primes in the
$p$-adic Lie extension, we compare the ranks of the $\pi$-free
quotient of the two dual Selmer groups. We then apply our results to compare the
characteristic elements attached to the Selmer groups. We also
study the variation of the ranks of the $\pi$-free quotient of the dual Selmer groups of specialization of a
big Galois representation. We emphasis that our results \textit{do
not} assume the vanishing of the $\mu$-invariant.
}
 \end{abstract}
 
 \medskip
\noindent Keywords and Phrases: Selmer groups, admissible
$p$-adic Lie extensions, $\M_H(G)$-property, $\pi$-free quotient.

\smallskip
\noindent Mathematics Subject Classification 2010: 11R23, 11R34, 11F80, 16S34.

\section{Introduction}
Throughout the paper, $p$ will always denote a rational prime. Let
$\Op$ be the ring of integers of a fixed finite extension $K$ of
$\Qp$, and fix a choice of a local parameter $\pi$ of $\Op$. Let $F$
be a number field and $F_{\infty}$ an admissible $p$-adic Lie
extension of $F$ with Galois group $G$. For a Galois representation
defined over a number field $F$ with coefficients in $\Op$, one can
attach a Selmer group to these data and this Selmer group carries a
natural module structure over the Iwasawa algebra $\Op\ps{G}$. The main
conjecture of Iwasawa theory is then a conjecture on a precise
relation between the Selmer group and a conjectural $p$-adic
$L$-function (for instance, see \cite{B,BV,CFKSV, EPW, FK, Gr89, Gr94, Kak,K,Oc, VBSD}).
In this paper, we are interested in comparing the Selmer groups of
two congruent Galois representations. One of the motivations behind
this study lies in the philosophy that the ``Iwasawa main
conjecture" should be preserved by congruences. Therefore, in view
of this philosophy, one expects that the various Iwasawa theoretical
invariants of the Selmer groups of two congruent Galois
representations should be related. Over the cyclotomic
$\Zp$-extension, such studies were carried out in \cite{AS, BS, EPW,
Gr94, GV, Ha, Oc, We}. Over a noncommutative $p$-adic Lie extension,
this has also been carried out in \cite{B,Ch09,CS12,Jh,LimAk,LimCMu, Sh,SS,SS14} to some extent.

In the context of a cyclotomic $\Zp$-extension, Greenberg and
Vatsal first considered the situation of two elliptic curves with
good ordinary reduction at the prime $p$, whose $p$-torsion points
are isomorphic as Galois modules. They proved that if the
$\mu$-invariant of one vanishes, so will the other, and they also
compared the $\la$-invariants of the Selmer groups (see
\cite[Theorem 1.4]{GV}). Following their footsteps, Emerton, Pollack
and Weston established similar results for Selmer groups of the
specializations of a Hida family (see \cite[Theorems 1 and 2]{GV}).
Since then, many authors have obtained similar results for more
general Galois representations (see \cite[Theorem 1.1]{Ha} and \cite[Theorems 1 and 2]{We}). In the situation of a noncommutative
$p$-adic Lie extension, analogous results have been obtained, and to
describe these results, we need to introduce some more terminology.
A Galois extension $F_{\infty}$ of $F$ is said to be an
\textit{admissible $p$-adic Lie extension} of $F$ if (i)
$G=\Gal(F_{\infty}/F)$ is compact $p$-adic Lie group, (ii)
$F_{\infty}$ contains the cyclotomic $\Zp$ extension $F^{\cyc}$ of
$F$ and (iii) $F_{\infty}$ is unramified outside a finite set of
primes. Write $H = \Gal(F_{\infty}/F^{\cyc})$. In \cite[Corollary
4.4.2]{B}, \cite[Theorems 4.2 and 4.11]{Ch09}, \cite[Proposition
15]{Jh}, \cite[Proposition 5.3]{LimAk} and \cite[Theorem 8.4]{SS},
it is proved that for two congruent Galois representations, if one
of the dual Selmer groups is finitely generated over $\Op\ps{H}$, so
is the other. Furthermore, they have also established relationship
between the $\Op\ps{H}$-ranks of two dual Selmer groups. A common
feature in all the above cited works is that it is shown that
whenever the Iwasawa $\mu$-invariant of one of the Selmer groups
vanishes, so does the other. It is then natural to consider the
situation when the said Iwasawa $\mu$-invariants are nonzero. In
this situation, it is convenient to divide the problem into two,
namely comparing the $\pi$-primary submodules of the dual Selmer
groups and comparing the $\pi$-free quotient of the dual Selmer
groups. For the comparison of the $\pi$-primary submodules of the
dual Selmer groups, the first attempt to make such a study was done
by Barman and A. Saikia (see \cite{BS}), where they compare the
$\mu$-invariants of the Selmer groups of two congruent elliptic
curves over the cyclotomic $\Zp$-extension. This result was later
extended to more general Galois representation over a noncommutative
admissible $p$-adic Lie extension in \cite{B, LimCMu}. It is
worthwhile mentioning that the results of the author in
\cite{LimCMu} proves a finer statement, namely, the $\pi$-primary
submodules of the dual Selmer groups of congruent Galois
representations have the same elementary representations (see
\cite[Theorem 4.2.1]{LimCMu}). For the comparison of the $\pi$-free
quotient of the dual Selmer groups, this was first considered by
Ahmed and Shekhar for two congruent elliptic curves over the
cyclotomic $\Zp$-extension (see \cite[Theorem 3.1]{AS}). In this
paper, we will consider the comparison of the $\pi$-free quotient of
the dual Selmer groups of an admissible $p$-adic Lie extension (and
for more general Galois representations).

We now like to mention another perspective of this paper. A classical conjecture of Mazur \cite{Mazur}
asserts that the dual Selmer group of an $p$-ordinary elliptic curve over the cyclotomic $\Zp$-extension is torsion. This conjecture is proven when the base field $F$ is abelian over $\Q$ (see \cite{K, R}). For an admissible $p$-adic Lie extension, the torsionness remains unknown except for some special situations (see \cite{HO,HV}). We also note that if the dual Selmer group is finitely
generated over $\Zp\ps{H}$, then it is automatically a torsion $\Zp\ps{G}$-module. Despite our lack of knowledge, we may still consider the following question. For two congruent Galois
representations, one may ask if one of the dual Selmer group is a
torsion $\Op\ps{G}$-module, is the other one also torsion? In the
context of \cite{AS, Ch09, EPW, GV, Jh}, this is not an issue, as
they already have the torsionness of Selmer groups by the results of
Rubin and Kato \cite{K, R}. In all previous works (for instance, \cite{Ha, LimAk,
SS, We}), the preservation of torsionness is established under the
stronger hypothesis that one of the dual Selmer group is a finitely
generated $\Op\ps{H}$-module. The preservation of the torsion property
under congruence (of a high enough power) has only been recently
obtained by the author in \cite[Theorem 4.2.1]{LimCMu}. (However, we emphasise that the property of dual Selmer group being torsion is far from being dependent on the
residual representation.) In this
paper, we like to go one step further, namely, we like to establish
the preservation of the so-called $\M_H(G)$-property under
congruence. In particular, in view of the discussion in the previous paragraph,
this paper can be thought as a complement to the results in \cite{LimCMu}, and at the same time, an extension of the results there.

A finitely generated $\Op\ps{G}$-module $M$ is said to satisfy the
\textit{$\M_H(G)$-property} if its $\pi$-free quotient
$M_f:=M/M(\pi)$ is finitely generated over $\Op\ps{H}$. Here
$M(\pi)$ denotes the $\pi$-primary submodule of $M$. It has been
conjectured for Galois representations coming from abelian varieties
with good ordinary reduction at $p$ or cuspidal eigenforms with good
ordinary reduction at $p$, the dual Selmer group associated to such
a Galois representation satisfies the $\M_H(G)$-property (see
\cite{CFKSV, CS12, Su}), and the validity of such a
property is necessary for the formulation of the main conjecture of
Iwasawa theory over a non-commutative $p$-adic Lie extension (see
\cite{B, BV, CFKSV, FK, Kak, VBSD}). At present, the only situation where the $\M_H(G)$ property
is known to hold is the ``$\mu =0$" situation (for instance, see \cite[Proposition 5.6]{CFKSV}, \cite[Theorem 2.1]{CS12} or \cite[Lemma 2.9]{Kak}). The verification of the $\M_H(G)$-property in general
seems out of reach at the moment (but see \cite[Section 2]{CSS}, \cite[Section 3]{CS12} or \cite[Section 3]{LimMHG} for some related discussion in this direction).  Granted this property and further
suppose that $H$ is pro-$p$ without $p$-torsion, it then makes sense
to speak of the $\Op\ps{H}$-rank of the $\pi$-free quotient of the
dual Selmer group, and this rank can be thought as a higher analog
of the classical $\la$-invariant in the cyclotomic situation (for
instance, see \cite{Ho}).

We now present our key result. Let $\big(A, \{A_v\}_{v|p},
\{A^+_v\}_{v|\R} \big)$ and $\big(B, \{B_v\}_{v|p},
\{B^+_v\}_{v|\R}\big)$ be two data which satisfy the conditions
(a)--(d) as in Section \ref{Arithmetic Preliminaries}. Then one can
attach dual Selmer groups to these data (see Section \ref{Arithmetic
Preliminaries}) which we denote by $X(A/F_{\infty})$ and
$X(B/F_{\infty})$. We now introduce the following congruence
condition on $A$ and $B$ which allows us to be able to compare the
Selmer groups of $A$ and $B$.

\smallskip \noindent $\mathbf{(Cong_n)}$ : There is an isomorphism
$A[\pi^{n}]\cong B[\pi^{n}]$ of $\Gal(\bar{F}/F)$-modules which induces a
$\Gal(\bar{F}_v/F_v)$-isomorphism $A_v[\pi^{n}]\cong B_v[\pi^{n}]$
for every $v|p$.

\smallskip
We define $e_G(A)$ to be
\[  \min\Big\{r~|~ \pi^r \big(X(A/F_{\infty})(\pi)\big)=0\Big\}. \]
Let $S$ denote a finite set of primes of $F$ which contains all the primes above
$p$, the ramified primes of $A$ and $B$, and the archimedean primes.
The following is our main theorem (see Theorem \ref{congruent thm} for a slightly refined version).

\begin{thm2} Let $F_{\infty}$ be a strongly admissible pro-$p$ $p$-adic Lie
extension of $F$ which is unramified outside $S$. Suppose that the following statements hold.
\begin{enumerate}
\item[$(1)$] $\mathbf{(Cong_{e_G(A)+1})}$ holds.
\item[$(2)$] $X(A/F_{\infty})$ satisfies the $\M_H(G)$-property.
\item[$(3)$] $X(B/F_{\infty})$ has no nonzero pseudo-null $\Op\ps{G}$-submodules.
   \item[$(4)$] For every $v\in S$, the decomposition group of $G$ at $v$ has dimension $\geq 2$.
   \end{enumerate}
    Then $X(B/F_{\infty})$ satisfies the $\M_H(G)$-property and we have
   \[ \rank_{\Op\ps{H}}\big(X_f(A/F_{\infty})\big) = \rank_{\Op\ps{H}}\big(X_f(B/F_{\infty})\big).\]
   \   \end{thm2}

The theorem can therefore be viewed as a generalization of the first
assertion of \cite[Theorem 4.2.1]{LimCMu}. We mention that the
equality of $\Op\ps{H}$-ranks in the theorem has been established in
\cite[Theorem 8.8]{SS} for every specialization of a big Galois
representation, but their proof relies on the much stronger
assumption that the dual Selmer group of the big Galois
representation satisfies the $\M_H(G)$-property (see Remark
\ref{SSremark}). We also note that our results apply to congruent
Galois representations that need not be specialization of a big
Galois representation, which their results do not apply.

We now say a little on the proof of the theorem, leaving the details
to the rest of the paper. In the special case when the admissible
$p$-extension is of dimension 2, the preservation of the
$\M_H(G)$-property has been established by the author in
\cite[Proposition 5.1.3]{LimCMu}. There we make use of a criterion
of Coates and Sujatha (cf.\ \cite[Corollary 3.2]{CS12}; also see
\cite[Lemma 5.1.1]{LimCMu}) which reduces the validity of
$\M_H(G)$-property to certain relations on the $\mu$-invariants. The
preservation of the $\M_H(G)$-property in this context then follows
by combining this criterion of Coates-Sujatha with the author's
comparison result (cf. \cite[Theorem 4.2.1]{LimCMu}) on the
structure of the $\pi$-primary submodule of the dual Selmer groups
of two congruent Galois representations over a general admissible
$p$-adic Lie extension. However, as the criterion of Coates and
Sujatha is only stated for an admissible $p$-extension of dimension
2, the above approach does not carry over for higher dimensional
$p$-adic Lie extensions. The obstruction towards extending the
criterion of Coates and Sujatha to a general admissible $p$-adic Lie
extension is that for a higher dimensional $p$-adic Lie extension,
one has to account for the higher $H$-homology groups of the dual
Selmer groups (for instance, see \cite[Theorem 3.1]{LimMHG}) which
vanishes in the situation of dimension 2. At our present knowledge,
it does not seem easy to study the structure of these homology
groups (but see \cite[Proposition 5.1]{LimMHG} or \cite[Proposition
4.3]{LimAk} for some discussion in this direction). Therefore, our
proof of the main theorem \textit{will} take a different route, and
the proof is, perhaps surprisingly, not difficult. Our approach is
inspired by \cite[Theorem 3.1]{AS} which compares the cyclotomic
$\la$-invariant of Selmer groups with positive $\mu$-invariant.
Following the said cited work, the idea is to first relate an
appropriate quotient of $X(A/F_{\infty})$ with the mod $\pi$
quotient of its $\pi$-free quotient (see Proposition \ref{alg
quotient compare}). Under the congruence condition, we then compare
the former quotients of $X(A/F_{\infty})$ and $X(B/F_{\infty})$ up
to finitely generated $k\ps{H}$-modules (here $k$ is the residue
field of $\Op$). Combining these with a Nakayama lemma argument,
this in turn allows us to deduce the validity of the
$\M_H(G)$-property of one from the other. The equality of the
$\Op\ps{H}$-ranks of the $\pi$-free quotient of $X(A/F_{\infty})$
and $X(B/F_{\infty})$ follows essentially the same argument with a
finer analysis.

We will apply our theorem to compare the Selmer groups of the
specializations of a big Galois representation. Namely, we prove
that for specializations in the big Galois representation which
satisfy appropriate congruence condition, the dual Selmer group of
one satisfies the $\M_H(G)$-property, then so will the other, and we
have an equality of the ranks of the $\pi$-free quotients of the
dual Selmer groups as predicted by Shekhar and Sujatha (see Theorem
\ref{big galois compare}).

In short, the results in this paper are concerned with comparing the
$\pi$-free quotient of the dual Selmer groups of congruent Galois
representations, and they complement the results proved by the
author in \cite{LimCMu} which is concerned with the comparison of
the $\pi$-primary submodules of the dual Selmer groups. The combined
results therefore give a rather satisfactory answer on comparing the
dual Selmer groups of two congruent Galois representations when the
congruence is a high power enough. At this point, it remains an open
problem whether one can do any meaningful comparison of the dual
Selmer groups when the Galois representations are congruent to each
other by a power lower than the exponent of the $\pi$-primary
submodules of the dual Selmer groups of the said Galois
representations. To the best of the author knowledge, this issue
does not seem to have been considered in literature, and
unfortunately, the author also does not have a mean to do this at
this point of writing. Another interesting problem is concerned with
the situation of a big Galois representation. Despite being able to
show the preservation of the $\M_H(G)$-property for appropriate
specializations of the big Galois representation, we are not able to
say anything on whether the dual Selmer group of the big Galois
representation satisfies the $\M_H(G)$-property (in the sense of
\cite[Section 5]{CS12}). In view of the noncommutative main
conjecture for big Galois representations in the sense of \cite{B},
we strongly believe that these two problems are important questions
of study and that the examination of these two problems will give
insights towards understanding the noncommutative main conjecture
for big Galois representations. We finally mention one more interesting problem
arising from our results. In view of our main theorem, it will be of interest to come up
with numerical examples. However, the obstructions to obtaining numerical results are that we do not have
a good understanding of the variation of $\mu$-invariants under base change and descent (see remark after Theorem \ref{congruent thm}). We hope to come back to these questions in subsequent papers.

We now give a brief description of the layout of the paper. In
Section \ref{Algebraic Preliminaries}, we recall certain algebraic
notion which will be used subsequently in the paper. It is here
where we develop a method to compare certain subquotients of two
$\Op\ps{G}$-modules which are annihilated by $\pi^{n+1}$ from some
$n$ (see Proposition \ref{alg compare}). We also identify certain
subquotient of a $\Op\ps{G}$-module with a quotient of its
$\pi$-free quotient module (see Proposition \ref{alg quotient
compare}), and this identification paves a way for us to apply a
Nakayama lemma argument. In Section \ref{Arithmetic Preliminaries},
we introduce the Selmer groups which are the main object of study in
this paper. Actually, to be precise, the Selmer group that we
consider is called the strict Selmer group in Greenberg's
terminology \cite{Gr89}. We also introduce another variant of
the Selmer group (called the Greenberg Selmer group) and an
appropriate Selmer complex which is closely related to the strict
Selmer group. In Section \ref{main results}, we will present and
prove our main results. Finally, Appendix \ref{pseudo-null
section} contains some discussion on the non-existence of
pseudo-null submodules of the dual (strict) Selmer group. Appendix
\ref{char subsection} contains some application of our main theorem
which compares the characteristic elements of $\pi$-free quotient of
the dual Selmer groups (see Theorem \ref{char congruent}) and this
can be thought as a refinement of a previous result of the author in
\cite[Theorem 6.3]{LimAk}. Although this latter result does not fit
into the theme of the paper, we have thought that it is interesting
enough to be included in an appendix.

\section{Algebraic Preliminaries} \label{Algebraic Preliminaries}

In this section, we establish some algebraic
preliminaries and notation which are necessary for us in order to
prepare for the discussion and the proofs of our results.

\subsection{Compact $p$-adic Lie group} \label{Compact p-adic Lie
group}

Fix a prime $p$. In this subsection, we recall some facts about
compact $p$-adic Lie groups. The standard references for the
material presented here are \cite{DSMS, Laz}.

For a finitely generated pro-$p$ group $G$, we write $G^{p^i} =
\langle g^{p^i}|~g\in G\rangle$, that is, the group generated by the
$p^i$th-powers of elements in $G$. The pro-$p$ group $G$ is said to
be \textit{powerful} if $G/\overline{G^{p}}$ is abelian for odd $p$,
or if $G/\overline{G^{4}}$ is abelian for $p=2$. If a powerful pro-$p$ group $G$
is torsionfree, we say that $G$ is \textit{uniform} (cf. \cite[Definition 4.1, Theorem 4.5]{DSMS}).

We now recall the following characterization of compact $p$-adic Lie
groups due to Lazard \cite{Laz} (see also \cite[Corollary 8.34]{DSMS}):
a topological group $G$ is a compact $p$-adic Lie group if and only
if $G$ contains a open normal uniform pro-$p$ subgroup. Furthermore,
if $G$ is a compact $p$-adic Lie group without $p$-torsion, it
follows from \cite[Corollaire 1]{Ser} (see also \cite[Chap.\ V
Sect.\ 2.2)]{Laz}) that $G$ has finite $p$-cohomological dimension.

\subsection{Torsion modules and pseudo-null modules} \label{Rank subsection}
As before, $p$ will denote a fixed prime.
Let $\Op$ be the ring of integers of a fixed finite extension of $\Qp$.
For a compact $p$-adic Lie group $G$, the
completed group algebra of $G$ over $\Op$ is given by
 \[ \Op\ps{G} = \plim_U \Op[G/U], \]
where $U$ runs over the open normal subgroups of $G$ and the inverse
limit is taken with respect to the canonical projection maps.

When $G$ is pro-$p$ and has no $p$-torsion, it is
well known that $\Op\ps{G}$ is an Auslander regular ring (cf.
\cite[Theorem 3.26]{V02} or \cite[Theorem A.1]{LimFine}), and has
no zero divisors (cf.\ \cite{Neu}). Therefore, $\Op\ps{G}$ admits a
skew field $K(G)$ which is flat over $\Op\ps{G}$ (see \cite[Chapters
6 and 10]{GW} or \cite[Chapter 4, \S 9 and \S 10]{Lam}). If $M$ is a
finitely generated $\Op\ps{G}$-module, we define the
$\Op\ps{G}$-rank of $M$ to be
$$ \rank_{\Op\ps{G}}(M)  = \dim_{K(G)} \big(K(G)\ot_{\Op\ps{G}}M\big). $$
The $\Op\ps{G}$-module $M$ is then said to be \textit{torsion} if
$\rank_{\Op\ps{G}} M = 0$. We will also make use of a well-known
equivalent definition for $M$ to be a torsion $\Op\ps{G}$-module,
namely: $\Hom_{\Op\ps{G}}(M, \Op\ps{G})=0$ (for instance, see
\cite[Lemma 4.2]{LimFine} or \cite[Lemma 2.2.1]{LimCMu}). A finitely
generated torsion $\Op\ps{G}$-module $M$ is said to be
\textit{pseudo-null} if $\Ext^1_{\Op\ps{G}}(M, \Op\ps{G}) =0$. For
an equivalent definition, we refer readers to \cite[Definitions 3.1
and 3.3; Proposition 3.5(ii)]{V02}. For the purpose of this article,
the definition we adopt will suffice. Finally, we mention that every
subquotient of a torsion $\Op\ps{G}$-module (resp., pseudo-null
$\Op\ps{G}$-module) is also torsion (resp. pseudo-null).

Now, fix a local parameter $\pi$ for $\Op$ and denote the residue
field of $\Op$ by $k$. The completed group algebra of $G$ over $k$
is given by
 \[ k\ps{G} = \plim_U k[G/U], \]
where $U$ runs over the open normal subgroups of $G$ and the inverse
limit is taken with respect to the canonical projection maps. For a
compact $p$-adic Lie group $G$ without $p$-torsion, it follows from
\cite[Theorem 3.30(ii)]{V02} (or \cite[Theorem A.1]{LimFine}) that
$k\ps{G}$ is an Auslander regular ring. Furthermore, if $G$ is
pro-$p$, then the ring $k\ps{G}$ has no zero divisors (cf.
\cite[Theorem C]{AB}). Therefore, one can define the notion of
$k\ps{G}$-rank as above when $G$ is pro-$p$ without $p$-torsion.
Similarly, we say that the module $N$ is a
\textit{torsion} $k\ps{G}$-module if $\rank_{k\ps{G}}N = 0$.

We end the subsection with some algebraic results which will be
used in the proof of our main result. As a start, we have the
following simple observation.

\bl \label{pi-module}
 Let $H$ be a compact $p$-adic Lie group and
$M$ a compact $\Op\ps{H}$-module. Suppose that $M$ is annihilated by
$\pi^{n+1}$. Then $M/M[\pi^n]$ and $M/\pi^n$ are $k\ps{H}$-modules.

Furthermore, if $M$ is finitely generated over $\Op\ps{H}$, then
$M/M[\pi^n]$ and $M/\pi^n$ are finitely generated over $k\ps{H}$.
\el

\bpf

It suffices to show that $M/M[\pi^n]$ and $M/\pi^n$ are annihilated
by $\pi$. Since $\pi^{n+1}M=0$, we have $\pi m\in M[\pi^{n}]$ for
all $m\in M$. This shows that $M/M[\pi^n]$ is annihilated by $\pi$.

By \cite[Lemma 5.2.5(b)]{RZ} and \cite[Proposition 7.4.1]{Wil}, we
have a surjection
  \[ \prod_{I}\Op\ps{H} \twoheadrightarrow M \]
  for some index set $I$. Since $M$ is annihilated by $\pi^{n+1}$, the surjection factors through a surjection
  \[ \prod_{I}\Op\ps{H}/\pi^{n+1} \twoheadrightarrow M.\]
It then follows from the following commutative diagram
\[ \SelectTips{eu}{}
\xymatrix{
  \prod_{I}\Op\ps{H}/\pi^{n+1}  \ar[d]^{\pi^n}
  \ar[r]^{} & M\ar[r] \ar[d]^{\pi^n} &0 \\
 \prod_{I}\Op\ps{H}/\pi^{n+1}\ar[r] & M
  \ar[r] &0   }
\]
that there is a surjection
 \[ \prod_{I}\Op\ps{H}/\pi \twoheadrightarrow M/\pi^n\]
 which in turn implies that $M/\pi^n$ is annihilated by $\pi$.

Finally, if $M$ is finitely generated over $\Op\ps{H}$, so are
$M/M[\pi^n]$ and $M/\pi^n$. But we have shown above that these two
modules are annihilated by $\pi$, and hence they are finitely
generated over $k\ps{H}$. \epf

For a closed subgroup $U$ of a compact $p$-adic Lie group $H$ and a
compact $\Op\ps{U}$-module $M$, we denote by $\Ind^U_H(M) :=
M\hat{\ot}_{\Op\ps{U}}\Op\ps{H}$ the compact induction of $M$ from
$U$ to $H$. If $M$ is a compact $k\ps{U}$-module, we can also view it as a compact
$\Op\ps{U}$-module and we have the identification $\Ind^U_H(M) =
M\hat{\ot}_{k\ps{U}}k\ps{H}$.

\bl \label{ind}
Let $H$ be a compact $p$-adic Lie group without $p$-torsion, $U$ a closed subgroup of $H$
and $M$ a finite $\Op\ps{U}$-module
which is annihilated by $\pi^{n+1}$. Then $\Ind^U_H(M)/\Ind^U_H(M)[\pi^{n}]$ and $\Ind^U_H(M)/\pi^n$ are finitely generated over $k\ps{H}$.

 Furthermore, if $U$ has dimension at least 1, then  $\Ind^U_H(M)/\big(\Ind^U_H(M)[\pi^{n}]\big)$ and $\Ind^U_H(M)/\pi^n$ are finitely generated
torsion $k\ps{H}$-modules.
\el

\bpf
 Clearly, $\Ind^U_H(M)$ is finitely generated over $\Op\ps{H}$ and is annihilated by $\pi^{n+1}$. Therefore, the first assertion follows from Lemma \ref{pi-module}. Now consider the following exact sequence
 \[ 0\lra M[\pi^n] \lra M\stackrel{\pi^n}{\lra} M \lra M/\pi^n\lra 0. \]
 Since $\Ind^U_H(-)$ is an exact functor (cf. \cite[Lemma 6.10.8]{RZ}), we have an exact sequence
 \[ 0\lra \Ind^U_H(M)[\pi^n] \lra \Ind^U_H(M)\stackrel{\pi^n}{\lra} \Ind^U_H(M) \lra \Ind^U_H(M)/\pi^n\lra 0.\]
It now follows easily from this exact sequence that
\[ \Ind^U_H(M)\big/\big(\Ind^U_H(M)[\pi^n]\big) \cong \Ind^U_H(M/M[\pi^n])\]
and \[ \Ind^U_H(M)/\pi^n \cong \Ind^U_H(M/\pi^n).\]
Since $M$ is finite, so are $M/M[\pi^n]$ and $M/\pi^n$. By Lemma \ref{pi-module}, they are finite $k\ps{U}$-modules. Now if $U$ has dimension at least 1, then $M/M[\pi^n]$ and $M/\pi^n$ are torsion $k\ps{U}$-modules. The second assertion then follows from an application of an $k$-analog of \cite[Lemma 5.5]{OcV02}.
\epf

The next proposition will be a key ingredient in our main theorem.

\bp \label{alg compare}
Let $H$ be a compact $p$-adic Lie group which has no $p$-torsion and has dimension at least 1. Suppose that we are given an exact sequence
\[ 0\lra C\lra M\lra N\lra D\lra 0\] of
$k\ps{H}$-modules which are annihilated by $\pi^{n+1}$, where $C$ is finitely generated over $\Op\ps{H}$ and $D$ is finite.
Then $M/M[\pi^n]$ is finitely generated over $k\ps{H}$ if and only if $N/N[\pi^n]$ is finitely generated over $k\ps{H}$.

Moreover, if $C/C[\pi^n]$ and $C/\pi^n$ are finitely generated torsion $k\ps{H}$-modules, then we have
\[ \rank_{k\ps{H}}M/M[\pi^n] = \rank_{k\ps{H}}N/N[\pi^n].\]
\ep

\bpf
Write $P=\im(M\lra N)$. We then have an exact sequence
\[ 0\lra P[\pi^n] \lra N[\pi^n] \lra D[\pi^n]\lra P/\pi^n. \]
Denoting $U = \im(N[\pi^n]\lra D[\pi^n])$, we then have the following commutative diagram
 \[ \SelectTips{eu}{}
\xymatrix{
  0 \ar[r] & P[\pi^n] \ar[d] \ar[r]^{} & N[\pi^n] \ar[d]
  \ar[r]^{} & U\ar[d] \ar[r]& 0 \\
  0 \ar[r] & P \ar[r] & N
  \ar[r] & D \ar[r] &0   }
\] with exact rows and the vertical maps are the inclusion maps. This in turn gives rise to the following short exact sequence
\[ 0\lra P/P[\pi^n] \lra N/N[\pi^n] \lra D/U\lra 0. \]
By Lemma \ref{pi-module}, this is an exact sequence of $k\ps{H}$-modules. Since $D$ is finite, so is $D/U$.
As $H$ has dimension $\geq 1$, it follows that $D/U$ is a finitely generated torsion $k\ps{H}$-module.
Therefore, we have that $P/P[\pi^n]$ is finitely generated over $k\ps{H}$ if and only if
$N/N[\pi^n]$ is finitely generated over $k\ps{H}$. Furthermore, in the event that these modules are finitely generated over $k\ps{H}$, we have
\[ \rank_{k\ps{H}}P/P[\pi^n] = \rank_{k\ps{H}}N/N[\pi^n].\]

Now consider the exact sequence
\[ 0\lra C[\pi^n]\lra M[\pi^n]\lra P[\pi^n] \lra C/\pi^n. \]
Denote $V= \im(M[\pi^n] \lra P[\pi^n]) $ and $W= \im(P[\pi^n]\lra C/\pi^n)$. We then have two commutative diagrams
\[ \SelectTips{eu}{}
\xymatrix{
  0 \ar[r] & C[\pi^n] \ar[d] \ar[r]^{} & M[\pi^n] \ar[d]
  \ar[r]^{} & V\ar[d] \ar[r]& 0 \\
  0 \ar[r] & C \ar[r] & M
  \ar[r] & P \ar[r] &0   }
\]
\[ \SelectTips{eu}{}
\xymatrix{
  0 \ar[r] & V \ar[d] \ar[r]^{} & P[\pi^n] \ar[d]
  \ar[r]^{} & W \ar[r]& 0 \\
  & P  \ar@{=}[r]  & P
   &  &  }
\]
with exact rows, where the vertical maps are the inclusion maps. From these two diagrams, we obtain two exact
sequences
\[ \ba{c}
 0\lra C/C[\pi^n] \lra M/M[\pi^n] \lra P/V \lra 0 \\
 0\lra W\lra P/V \lra P/P[\pi^n]\lra 0.
 \ea \]
Since $C$
is finitely generated over $\Op\ps{H}$, it follows from Lemma \ref{pi-module}
that $ C/\pi^n$ and $C/C[\pi^n]$ are finitely generated over $k\ps{H}$. As $W\subseteq C/\pi^n$, $W$ is also finitely generated over $k\ps{H}$.
Hence it follows from the two exact sequences that $M/M[\pi^n]$ is
finitely generated over $k\ps{H}$ if and only if $P/P[\pi^n]$ is
finitely generated over $k\ps{H}$. Combining this observation with
that in the previous paragraph, we have that $M/M[\pi^n]$ is
finitely generated over $k\ps{H}$ if and only if $N/N[\pi^n]$ is
finitely generated over $k\ps{H}$.

Now if we assume further that $C/C[\pi^n]$ and $C/\pi^n$ are
finitely generated torsion $k\ps{H}$-modules, it then follows that
$W$ is also a finitely generated torsion $k\ps{H}$-module. By
analysing the above two exact sequences again, we have
\[ \rank_{k\ps{H}}M/M[\pi^n] = \rank_{k\ps{H}}P/P[\pi^n].\]
Combining this with the equality obtained at the end of the first
paragraph of the proof, we obtained the desired equality. \epf

We record another useful lemma.

\bl \label{useful lemma} Let $H$ be a compact $p$-adic Lie group
which has no $p$-torsion and has dimension at least 1. Let $M$ be a
finitely generated $\Op\ps{H}$-module which is annihilated by
$\pi^{n+1}$. Suppose that
$M/M[\pi^n]$ and $M/\pi^n$ are torsion $k\ps{H}$-modules. Then for
every $\Op\ps{H}$-subquotient $N$ of $M$, $N/N[\pi^n]$ and $N/\pi^n$
are torsion $k\ps{H}$-modules. \el

\bpf It suffices to prove the assertion for the cases when $N$ is a
submodule of $M$ and when $N$ is a quotient of $M$. We first suppose that
$M\tha N$. It is then straightforward to verify that this surjection induces surjections $ M/M[\pi^n]\tha N/N[\pi^n]$ and
$M/\pi^n\tha N/\pi^n$. Hence $N/N[\pi^n]$ and $N/\pi^n$
are torsion $k\ps{H}$-modules.

Now suppose that $N\subseteq M$. Then one can check easily that
$N/N[\pi^n]\subseteq M/M[\pi^n]$, and so $N/N[\pi^n]$ is a torsion
$k\ps{H}$-module. It therefore remains to show that $N/\pi^n$ is a
torsion $k\ps{H}$-module. Consider the following exact sequence
\[ 0\lra N[\pi^n]\lra M[\pi^n]\lra P[\pi^n] \lra N/\pi^n \lra M/\pi^n\lra P/\pi^n\lra 0. \]
Denote $V= \im(M[\pi^n] \lra P[\pi^n]) $ and $W= \im(P[\pi^n]\lra
N/\pi^n)$. By a similar argument to that in Proposition \ref{alg
compare}, we obtain three exact sequences
\[ \ba{c}
   0\lra N/N[\pi^n]\lra M/M[\pi^n] \lra P/V\lra 0, \\
   0\lra W\lra P/V\lra P/P[\pi^n] \lra 0, \\
   0\lra W\lra N/\pi^n \lra M/\pi^n \lra P/\pi^n \lra 0
   \ea\]
of $k\ps{H}$-modules. (Note that by Lemma \ref{pi-module},
$M/M[\pi^n]$ is a $k\ps{H}$-module, and hence, so are $P/V$ and
$W$.) Since $M/M[\pi^n]$ is a torsion $k\ps{H}$-module, it follows
from the first two exact sequences that $W$ is a torsion
$k\ps{H}$-module. Combining this with the hypothesis that $M/\pi^n$
is a torsion $k\ps{H}$-module, it follows from the third exact
sequence that $N/\pi^n$ is a torsion $k\ps{H}$-module. \epf

\subsection{$\mu$-invariant and $\pi$-primary modules} \label{pi-primary modules}

For a given finitely generated $\Op\ps{G}$-module $M$, we denote by
$M(\pi)$ the $\Op\ps{G}$-submodule of $M$ which consists of
elements of $M$ that are annihilated by some power of $\pi$. Since
the ring $\Op\ps{G}$ is Noetherian, the module $M(\pi)$ is certainly
finitely generated over $\Op\ps{G}$. Hence one can find an integer
$r\geq 0$ such that $\pi^r$ annihilates $M(\pi)$. The
\textit{$\pi$-exponent} of $M$ is then defined to be
\[ e_{\Op\ps{G}}(M) = \min\{r~|~ \pi^rM(\pi)=0\}. \]
Now suppose that $G$ is pro-$p$ without $p$-torsion. Following \cite[Formula (33)]{Ho}, we define the
\textit{$\mu$-invariant}
  \[\mu_{\Op\ps{G}}(M) = \sum_{i\geq 0}\rank_{k\ps{G}}\big(\pi^i
   M(\pi)/\pi^{i+1}\big). \]
(For another alternative, but equivalent, definition, see
\cite[Definition 3.32]{V02}.) By the above discussion and our
definition of $k\ps{G}$-rank, the sum on the right is a finite one.
It is clear from the definition that $\mu_{\Op\ps{G}}(M) =
\mu_{\Op\ps{G}}(M(\pi))$. Also, it is not difficult to see that this
definition coincides with the classical notion of the
$\mu$-invariant for $\Ga$-modules when $G=\Ga\cong\Zp$.

Continue supposing that $G$ is pro-$p$ without $p$-torsion. Then both $\Op\ps{G}$ and $k\ps{G}$ are Auslander
regular rings with no zero divisors. For a finitely generated
$\Op\ps{G}$-module $M$, it then follows from \cite[Proposition
1.11]{Ho2} (see also \cite[Theorem 3.40]{V02}) that there is a
$\Op\ps{G}$-homomorphism
\[ \varphi: M(\pi) \lra \bigoplus_{i=1}^s\Op\ps{G}/\pi^{\al_i},\] whose
kernel and cokernel are pseudo-null $\Op\ps{G}$-modules, and where
the integers $s$ and $\al_i$ are uniquely determined. We will call
$\bigoplus_{i=1}^s\Op\ps{G}/\pi^{\al_i}$  the \textit{elementary
representation} of $M(\pi)$. In fact, in the process of establishing
the above, one also has the equality $\mu_{\Op\ps{G}}(M) = \displaystyle
\sum_{i=1}^s\al_i$ (see loc. cit.). We set
\[\theta_{\Op\ps{G}}(M) : = \max_{1\leq i\leq s}\{\al_i\}.\]
It is not difficult to see that $e_{\Op\ps{G}}(M) \geq \theta_{\Op\ps{G}}(M)$. The following lemma gives a sufficient criterion for equality to hold.

\bl \label{e=theta} Let $G$ be a compact $p$-adic Lie group which is
pro-$p$ and has no $p$-torsion. Suppose that $M$ is a finitely
generated $\Op\ps{G}$-module which has no nonzero pseudo-null
$\Op\ps{G}$-submodules. Then
\[ e_{\Op\ps{G}}(M) = \theta_{\Op\ps{G}}(M).\] \el

\bpf
Let
\[ \varphi: M(\pi) \lra \bigoplus_{i=1}^s\Op\ps{G}/\pi^{\al_i}\]
be an $\Op\ps{G}$-homomorphism, whose kernel and cokernel are pseudo-null $\Op\ps{G}$-modules.
Since $M$ has no nonzero pseudo-null $\Op\ps{G}$-submodules, it follows that $\varphi$ is injective, and hence
$\pi^{\theta_{\Op\ps{G}}(M)}$ annihilates $M(\pi)$. This in turn yields the required equality.
\epf

Finally, we introduce the notion of the $\M_H(G)$-property. Let $G$
be a compact $p$-adic Lie group with a closed subgroup $H$ such that
$G/H\cong \Zp$. We then say that a finitely generated
$\Op\ps{G}$-module $M$ \textit{satisfies the $\M_H(G)$-property} if
its $\pi$-free quotient $M_f:=M/M(\pi)$ is finitely generated over
$\Op\ps{H}$. As noted in the introductional section, it has been
conjectured for certain Galois representations coming from abelian
varieties with good ordinary reduction at $p$ or cuspidal eigenforms
with good ordinary reduction at $p$, the dual Selmer group
associated to such a Galois representation satisfies the
$\M_H(G)$-property (see \cite{CFKSV, CS12, FK, Su}).

We end with another important proposition which will play a part in
the proof of our main theorem. As mentioned in the introduction, the
proposition and its proof are inspired by the proof of \cite[Theorem
3.1]{AS}. The point of the proposition is to relate an appropriate
quotient of $M$ to the mod-$\pi$ quotient of its $\pi$-free
quotient.

\bp \label{alg quotient compare}
Let $G$ be a compact $p$-adic Lie group
and $H$ a closed subgroup of $G$ such that $G/H\cong \Zp$.
 Let $M$ be a finitely generated torsion $\Op\ps{G}$-module.
  Then for every $n\geq 1$, we have a short exact sequence
  \[ 0\lra \frac{M(\pi)/\pi^{n+1}}{(M(\pi)/\pi^{n+1})[\pi^n]} \lra \frac{M/\pi^{n+1}}{(M/\pi^{n+1})[\pi^n]} \lra M_f/\pi \lra 0\]
  of finitely generated $k\ps{G}$-modules (and hence compact $k\ps{H}$-modules).
  In particular, when $n\geq e_{\Op\ps{G}}(M)$, we have an isomorphism
 \[ \frac{M/\pi^{n+1}}{(M/\pi^{n+1})[\pi^n]} \cong M_f/\pi\]
 of finitely generated $k\ps{G}$-modules (and hence compact $k\ps{H}$-modules).
\ep

\bpf
Consider the following commutative
  diagram
  \[ \SelectTips{eu}{}
\xymatrix{
  0 \ar[r] & M(\pi) \ar[d]^{\pi^n} \ar[r]^{} & M \ar[d]^{\pi^n}
  \ar[r]^{} & M_f \ar[d]^{\pi^n} \ar[r]& 0 \\
  0 \ar[r] & M(\pi) \ar[r] & M
  \ar[r] & M_f \ar[r] &0   }
\] with exact rows, where the vertical maps are given by
multiplication by $\pi^n$. Since $M_f$ has no $\pi$-torsion, the
rightmost vertical map is injective. Therefore, it follows that there is a short exact sequence
\[ 0\lra M(\pi)/\pi^n \lra M/\pi^n \lra M_f/\pi^n \lra 0\]
of $\Op\ps{G}$-modules. We also have a similar short exact sequence replacing $n$ by $n+1$. This latter short exact sequence gives rise to a long exact sequence
\[ 0\lra  \big(M(\pi)/\pi^{n+1}\big)[\pi^n] \lra \big(M/\pi^{n+1}\big)[\pi^n] \lra \big(M_f/\pi^{n+1}\big)[\pi^n] \hspace{1in}\]
\[\hspace{1.5in}\lra  M(\pi)/\pi^n \lra M/\pi^n \lra M_f/\pi^n\lra  0. \]
Since the last three terms of the long exact sequence is part of the former short exact sequence, we deduce that the first three terms of the long exact sequence actually form a short exact sequence
\[  0\lra  \big(M(\pi)/\pi^{n+1}\big)[\pi^n] \lra \big(M/\pi^{n+1}\big)[\pi^n] \lra \big(M_f/\pi^{n+1}\big)[\pi^n] \lra 0.\]
From this, we have a short exact sequence
\[ 0\lra \frac{M(\pi)/\pi^{n+1}}{(M(\pi)/\pi^{n+1})[\pi^n]} \lra \frac{M/\pi^{n+1}}{(M/\pi^{n+1})[\pi^n]} \lra \frac{M_f/\pi^{n+1}}{(M_f/\pi^{n+1})[\pi^n]}\lra 0 \]
of $k\ps{G}$-modules (noting Lemma \ref{pi-module}). It remains to show that
\[\frac{M_f/\pi^{n+1}}{(M_f/\pi^{n+1})[\pi^n]}\cong M_f/\pi. \]
Indeed, as $M_f$ has no $\pi$-torsion, one can easily check that
\[(M_f/\pi^{n+1})[\pi^n] = \pi M_f/\pi^{n+1},\]
and the required isomorphism is immediate from this.

Finally, if $n\geq e_{\Op\ps{G}}(M)$, then
\[\frac{M(\pi)/\pi^{n+1}}{(M(\pi)/\pi^{n+1})[\pi^n]} = \frac{M(\pi)}{M(\pi)} =0.\]
The proof of the proposition is now completed.
\epf

We record one more result which gives an upper bound of the $\pi$-exponent.

\bp \label{upper bound}
Let $G$ be a compact $p$-adic Lie group
and $H$ a closed subgroup of $G$ such that $G/H\cong \Zp$.
Let $M$ be a finitely generated torsion $\Op\ps{G}$-module which has no nonzero pseudo-null
$\Op\ps{G}$-submodules. Suppose that $\frac{M/\pi^{n+1}}{(M/\pi^{n+1})[\pi^n]}$ is finitely generated over $k\ps{H}$ for some $n$. Then we have $\theta_{\Op\ps{G}}(M)=e_{\Op\ps{G}}(M)\leq n$. In particular, we have an isomorphism
 \[ \frac{M/\pi^{n+1}}{(M/\pi^{n+1})[\pi^n]} \cong M_f/\pi\]
 of finitely generated $k\ps{G}$-modules $($and hence compact $k\ps{H}$-modules$)$.
 \ep

\bpf The first equality is a consequence of Lemma \ref{e=theta}. We now proceed with proving $\theta_{\Op\ps{G}}(M)\leq n$. By Proposition \ref{alg quotient compare},
$\frac{M(\pi)/\pi^{n+1}}{(M(\pi)/\pi^{n+1})[\pi^n]}$ is finitely
generated over $k\ps{H}$, and hence over $\Op\ps{H}$. Thus, this
module has trivial $\mu_{\Op\ps{G}}$-invariant (cf.\ \cite[Lemma
2.7]{Ho}). This in turn implies
that
\[ \mu_{\Op\ps{G}} \Big(M(\pi)/\pi^{n+1}\Big)  = \mu_{\Op\ps{G}}\Big( \big(M(\pi)/\pi^{n+1}\big)[\pi^n] \Big).\]
Consider an $\Op\ps{G}$-homomorphism
\[ \varphi: M(\pi) \lra \bigoplus_{i=1}^s\Op\ps{G}/\pi^{\al_i},\] whose
kernel and cokernel are pseudo-null $\Op\ps{G}$-modules. Then it is
straightforward to verify that $\varphi$ induces two
$\Op\ps{G}$-homomorphisms
\[ M(\pi)/\pi^{n+1} \lra \bigoplus_{i=1}^s\Op\ps{G}/\pi^{\min\{n+1,\al_i\}}\]
and
\[ \big(M(\pi)/\pi^{n+1}\big)[\pi^n] \lra \bigoplus_{i=1}^s\Op\ps{G}/\pi^{\min\{n,\al_i\}},\]
whose kernels and cokernels are pseudo-null $\Op\ps{G}$-modules.
Since the two modules in question have the same
$\mu_{\Op\ps{G}}$-invariants by the above discussion, we have an
equality
\[\sum_{i=1}^s\min\{n,\al_i\} =\sum_{i=1}^s\min\{n+1,\al_i\}\]
which in turn implies that $\al_i\leq n$ for all $i$. Hence
$e_{\Op\ps{G}}(M)=\theta_{\Op\ps{G}}(M)\leq n$. This proves the first assertion of the proposition.
The second assertion of the proposition is
then immediate from this and Proposition \ref{alg quotient
compare}.\epf

\section{Arithmetic Preliminaries} \label{Arithmetic Preliminaries}

In this section, we introduce the Selmer groups and Selmer
complexes. At the same time, we fix the notation that we shall use throughout
the paper.

\subsection{Arithmetc datum}
To start, let $p$ be a prime and $F$ a number
field. If $p=2$, assume further that our number field $F$ has no real primes.
Denote by $\Op$ the ring of integers of some finite extension $K$
of $\Qp$. We then fix a local parameter $\pi$ for $\Op$. Suppose that we
are given the following datum  $\big(A, \{A_v\}_{v|p},
\{A^+_v\}_{v|\R} \big)$ defined over $F$:

\begin{enumerate}
 \item[(a)] $A$ is a
cofree $\Op$-module of $\Op$-corank $d$ with a
continuous, $\Op$-linear $\Gal(\bar{F}/F)$-action which is
unramified outside a finite set of primes of $F$.

 \item[(b)] For each prime $v$ of $F$ above $p$, $A_v$ is a
$\Gal(\bar{F}_v/F_v)$-submodule of $A$ which is cofree of
$\Op$-corank $d_v$.

\item[(c)] For each real prime $v$ of $F$, we write $A_v^+=
A^{\Gal(\bar{F}_v/F_v)}$  which is assumed to be cofree of
$\Op$-corank $d^+_v$.

\item[(d)] The quantities $d, d_v$ and $d_v^+$ satisfy the following identity
  \begin{equation} \label{data equality}
  \sum_{v|p} (d-d_v)[F_v:\Qp] = dr_2(F) +
 \sum_{v~\mathrm{real}}(d-d^+_v),
  \end{equation}
  where $r_2(F)$ is the number of complex primes of $F$.
\end{enumerate}

 We now consider the base change property of our datum. Let $L$ be a
finite extension of $F$. We can then obtain another datum $\big(A,
\{A_w\}_{w|p}, \{A^+_w\}_{w|\R} \big)$ over $L$ as follows: we
consider $A$ as a $\Gal(\bar{F}/L)$-module, and for each prime $w$
of $L$ above $p$, we set $A_w =A_v$, where $v$ is a prime of $F$
below $w$, and view it as a $\Gal(\bar{F}_v/L_w)$-module. Then $d_w
= d_v$. For each real prime $w$ of $L$, one sets
$A^{\Gal(\bar{L}_w/L_w)}= A^{\Gal(\bar{F}_v/F_v)}$ and writes $d^+_w
= d^+_v$, where $v$ is a real prime of $F$ below $w$. In general,
the $d_w$'s and $d_w^+$'s need not satisfy equality (\ref{data
equality}). We now record the following lemma which gives some
sufficient conditions for equality (\ref{data equality}) to hold for
the datum $\big(A, \{A_w\}_{w|p}, \{A^+_w\}_{w|\R} \big)$ over $L$.

\bl \label{data base change} Suppose that $\big(A, \{A_v\}_{v|p},
\{A^+_v\}_{v|\R} \big)$ is a datum defined over $F$.
 Suppose further that at least one of the following statements holds.
 \begin{enumerate}
\item[$(i)$] All the archimedean primes of $F$ are unramified in
$L$.

\item[$(ii)$] $[L:F]$ is odd

\item[$(iii)$] $F$ is totally imaginary.

\item[$(iv)$] $F$ is totally real, $L$ is totally imaginary and
\[ \sum_{v~\mathrm{real}} d^+_v = d[F:\Q]/2.\]
 \end{enumerate}
Then we have the equality
 \[ \sum_{w|p} (d-d_w)[L_w:\Qp] = dr_2(L) +
 \sum_{w~\mathrm{real}}(d-d^+_w).\]
 \el

\bpf Since the only ramified archimedean primes are real primes and real
primes can only ramify in an extension of even degree, it follows
that if either of the assertions in (ii) or (iii) holds, then the
assertion in (i) holds. Therefore, to prove the lemma in these
cases, it suffices to prove it under the assumption of (i). We first
perform the following calculation
 \[ \ba{rl}
\displaystyle \sum_{w|p} (d-d_w)[L_w:\Qp]\!\!
   &= \displaystyle\sum_{v|p}\sum_{w|v} (d-d_v)[L_w:F_v][F_v:\Qp] \\
   &= \displaystyle\sum_{v|p} (d-d_v)[F_v:\Qp] \sum_{w|v}[L_w:F_v] \\
   &= \displaystyle [L:F]\sum_{v|p} (d-d_v)[F_v:\Qp] \\
  &= \displaystyle [L:F]\Big(dr_2(F) + \sum_{v~\mathrm{real}} (d-d^+_v)\Big) \\
  &= \displaystyle d[L:F]r_2(F) + [L:F]\sum_{v~\mathrm{real}} (d-d^+_v). \\  \ea \]
 Now if (i) holds, then every prime of $L$ above a real prime (resp.,
complex prime) of $F$ is a real prime (resp., complex prime).
Therefore, one has $[L:F]r_2(F) = r_2(L)$ and
\[ [L:F]\sum_{v~\mathrm{real}} (d-d^+_v) = \sum_{w~\mathrm{real}} (d-d^+_w). \]
The required conclusion then follows.

Now suppose that (iv) holds. Then $r_2(F) =0$ and we have
 \[ \ba{rl}
\displaystyle \sum_{w|p} (d-d_w)[L_w:\Qp] =\displaystyle
[L:F]\sum_{v~\mathrm{real}} (d-d^+_v)\!\! &=
  \displaystyle[L:F]\sum_{v~\mathrm{real}}d - [L:F]\sum_{v~\mathrm{real}}d^+_v \\
   &= [L:F][F:\Q]d - [L:F] d[F:\Q]/2 \vspace{0.1in} \\
 &= d[L:\Q]/2 = dr_2(L).
\ea
\]
\epf

We now describe the arithmetic situation, where we can
obtain the datum from. Let $V$ be a $d$-dimensional
$K$-vector space with a continuous $\Gal(\bar{F}/F)$-action which is
unramified outside a finite set of primes. Suppose that for each
prime $v$ of $F$ above $p$, there is a $d_v$-dimensional
$K$-subspace $V_v$ of $V$ which is invariant under the action of
$\Gal(\bar{F}_v/F_v)$, and for each real prime $v$ of $F$,
$V^{\Gal(\bar{F}_v/F_v)}$ has dimension $d_v^+$. Choose a
$\Gal(\bar{F}/F)$-stable $\Op$-lattice $T$ of $V$ (such a lattice
exists by compactness). We can obtain a data as above from $V$ by
setting $A= V/T$ and $A_v = V_v/ (T\cap V_v)$. Note that both $A$ and $A_v$
depend on the choice of the lattice $T$. The basic examples of such
Galois representations are (1) $V=V_p(E)$, where $E$ is an elliptic
curve with either good ordinary reduction or multiplicative
reduction at each prime of $F$ above $p$, and (2) $V$ is the Galois
representation attached to a primitive Hecke eigenform which is
ordinary at $p$. For more examples of how the above datum arises
from Galois representations, we refer readers to \cite[\S 9]{Gr89},
\cite[Section 1.2]{We} and \cite[Section 3]{LimCMu}.

\subsection{Selmer groups}

We now introduce two variants of Selmer groups due to Greenberg \cite{Gr89}. Let $S$ be a finite set of
primes of $F$ which contains all the primes above $p$, the ramified
primes of $A$ and all the infinite primes of $F$. Denote by $F_S$ the
maximal algebraic extension of $F$ unramified outside $S$ and write
$G_S(\mathcal{L}) = \Gal(F_S/\mathcal{L})$ for every algebraic
extension $\mathcal{L}$ of $F$ which is contained in $F_S$. Let $L$
be a finite extension of $F$ contained in $F_S$ such that the datum
$\big(A, \{A_w\}_{w|p}, \{A^+_w\}_{w|\R} \big)$ satisfies (\ref{data
equality}). For a prime $w$ of $L$ lying over $S$, set
\[ H^1_{str}(L_w, A)=
\begin{cases} \ker\big(H^1(L_w, A)\lra H^1(L_w, A/A_w)\big) & \text{\mbox{if} $w$
 divides $p$},\\
 \ker\big(H^1(L_w, A)\lra H^1(L^{ur}_w, A)\big) & \text{\mbox{if} $w$ does not divide $p$,}
\end{cases} \]
 where $L_w^{ur}$ is the maximal unramified extension of $L_w$.
The (strict) Selmer group attached to the datum is then defined by
\[ S(A/L) := \Sel^{str}(A/L) := \ker\Big( H^1(G_S(L),A)\lra
\bigoplus_{w \in S_L}H^1_s(L_w, A)\Big),\] where we write
$H^1_s(L_w, A) = H^1(L_w, A)/H^1_{str}(L_w, A)$ and $S_L$ denotes
the set of primes of $L$ above $S$. We then write $X(A/L)$ for the
Pontryagin dual of $S(A/L)$.

\smallskip
 A Galois extension $F_{\infty}$ of $F$ is said to be an
\textit{$S$-admissible $p$-adic Lie extension} of $F$ if (i)
$\Gal(F_{\infty}/F)$ is compact $p$-adic Lie group, (ii)
$F_{\infty}$ contains the cyclotomic $\Zp$ extension $F^{\cyc}$ of
$F$ and (iii) $F_{\infty}$ is contained in $F_S$. Write $G =
\Gal(F_{\infty}/F)$, $H = \Gal(F_{\infty}/F^{\cyc})$ and $\Ga
=\Gal(F^{\cyc}/F)$. In the event that $\Gal(F_{\infty}/F)$ is a compact
$p$-adic Lie group without $p$-torsion, we say that $F_{\infty}$ is a \textit{strongly
$S$-admissible $p$-adic Lie extension} of $F$.

We define $S(A/F_{\infty}) = \ilim_L S(A/L)$,
where the limit runs over all finite extensions $L$ of $F$ contained
in $F_{\infty}$. We shall write $X(A/F_{\infty})$ for the Pontryagin
dual of $S(A/F_{\infty})$. By a similar argument to that in
\cite[Corollary 2.3]{CS12}, one can show that $X(A/F_{\infty})$ is
independent of the choice of $S$ as long as $S$ contains all the
primes above $p$, the ramified primes of $A$, the primes that ramify
in $F_{\infty}/F$ and all infinite primes.

We introduce another variant of the Selmer group which is usually
called the Greenberg Selmer group. Set
\[ H^1_{Gr}(F_v, A)=
\begin{cases} \ker\big(H^1(F_v, A)\lra H^1(F_v^{ur}, A/A_v)\big) & \text{\mbox{if} $v|p$},\\
 \ker\big(H^1(F_v, A)\lra H^1(F^{ur}_v, A)\big) & \text{\mbox{if} $v\nmid p$.}
\end{cases} \]
The Greenberg Selmer group attached to the datum $\big(A,
\{A_w\}_{w|p}, \{A^+_w\}_{w|\R} \big)$  is then defined by
\[ \Sel^{Gr}(A/F) = \ker\Big( H^1(G_S(F),A)\lra \bigoplus_{v \in S}H^1_g(F_v,
A)\Big),\] where we write $H^1_g(F_v, A) = H^1(F_v, A)/H^1_{Gr}(F_v,
A)$. For an $S$-admissible $p$-adic Lie extension $F_{\infty}$, we
define $\Sel^{Gr}(A/F_{\infty}) = \ilim_L \Sel^{Gr}(A/L)$ and denote by
$X^{Gr}(A/F_{\infty})$ the Pontryagin dual of
$\Sel^{Gr}(A/F_{\infty})$.

\medskip
We end the subsection by comparing the two
Selmer groups of Greenberg which will take the form of two lemmas.

\bl \label{Greenberg compare lemma}
  We have an exact sequence
\[0\lra N \lra X^{Gr}(A/F_{\infty}) \lra X(A/F_{\infty})\lra 0,\]
where $N$ is a finitely generated $\Op\ps{H}$-module.

Suppose further that for each $v|p$, the
decomposition group of $G=\Gal(F_{\infty}/F)$ at $v$ has dimension
$\geq 2$. Then $N$ is a finitely generated torsion $\Op\ps{H}$-module. \el

\bpf
Consider the following commutative
diagram
\[  \entrymodifiers={!! <0pt, .8ex>+} \SelectTips{eu}{}\xymatrix{
    0 \ar[r]^{} & S(A/F_{\infty}) \ar[d] \ar[r] &  H^1(G_S(F_{\infty}), A) \ar@{=}[d]
    \ar[r] & \bigoplus_{v \in S} J_v(A/F_{\infty}) \ar[d]^{\al} \\
    0 \ar[r]^{} & \Sel^{Gr}(A/F_{\infty}) \ar[r]^{} & H^1(G_S(F_{\infty}), A) \ar[r] & \
    \bigoplus_{v \in S} J^{Gr}_v(A/F_{\infty})  } \]
with exact rows, where $J^{Gr}_v(A/F_{\infty}) =
\ilim_L\bigoplus_{w|v} H^1_g(L_w, A)$. We first show
that $\ker \al$ is cofinitely generated over $\Op\ps{H}$. Write $\al
=\oplus_v \al_v$, where $v$ runs over the set of primes of $F$
above $S$. Clearly,
$J_v(A/F_{\infty}) = J^{Gr}_v(A/F_{\infty})$ for $v\nmid p$. Therefore, $\ker\al_v =0$ for these $v$'s.

Now for each $v|p$, write $\al_v = \oplus\al_w$, where $w$ runs over the set of primes of $F^{\cyc}$ above $v$.
Fix a prime of $F_{\infty}$ above $v$ which we also denote by $w$. Write
$I_{\infty, w}$ for the inertia subgroup of
$\Gal(\overline{F}_{\infty,w}/F_{\infty, w})$ and $U_w =
\Gal(\overline{F}_{\infty,w}/F_{\infty, w})/I_{\infty, w}$. It then
follows from the Hochschild-Serre spectral sequence that we have
\[ 0\lra H^1(U_w, (A/A_v)^{I_{\infty,w}}) \lra H^1(F_{\infty, w}, A/A_v)
\lra H^1(I_{\infty,w}, A/A_v)^{U_w}.\]
 Since $U_w$ is topologically cyclic, $H^1(U_w,
(A/A_v)^{I_{\infty,w}}) \cong
\big((A/A_v)^{I_{\infty,w}}\big)_{U_v}$ and so is cofinitely
generated over $\Op$. Let $H_w$ denote the decomposition subgroup of
$H$ corresponding to $w$. Then $\ker\al_w =
\mathrm{Coind}^{H_w}_H\big(H^1(U_w, (A/A_v)^{I_{\infty,w}})\big)$ is
cofinitely generated over $\Op\ps{H}$. Since there is only finite
number of primes of $F^{\cyc}$ above $v$, we have that $\ker\al_v$
is cofinitely generated over $\Op\ps{H}$.

Now suppose that for each $v|p$, the decomposition group of
$G=\Gal(F_{\infty}/F)$ at $v$ has dimension $\geq 2$. Then $H_w$ has
dimension $\geq 1$, and it follows from this that $H^1(U_w,
(A/A_v)^{I_{\infty,w}})$ is a cofinitely generated torsion
$\Op\ps{H_w}$-module. By \cite[Lemma 5.5]{OcV02}, this in turn
implies that $\ker\al_v$ is a cofinitely generated torsion
$\Op\ps{H}$-module. \epf

In the next lemma, we will write $X_f(A/F_{\infty})$ (resp.,
$X_f^{Gr}(A/F_{\infty})$) for the $\pi$-free quotient of
$X(A/F_{\infty})$ (resp., the $\pi$-free quotient of
$X^{Gr}(A/F_{\infty})$).

\bl \label{Greenberg equal}
 $X(A/F_{\infty})$ satisfies the $\M_H(G)$-property if and only if $X^{Gr}(A/F_{\infty})$
satisfies the $\M_H(G)$-property.

Suppose further that $G=\Gal(F_{\infty}/F)$ is pro-$p$ without $p$-torsion and that for each $v|p$, the
decomposition group of $G$ at $v$ has dimension
$\geq 2$. Then if $X(A/F_{\infty})$ $($and hence $X^{Gr}(A/F_{\infty}))$ satisfies the $\M_H(G)$-property,
we have
\[ \rank_{\Op\ps{H}}\big(X_f(A/F_{\infty})\big) = \rank_{\Op\ps{H}}\big(X^{Gr}_f(A/F_{\infty})\big).\]
 \el

\bpf We begin by verifying the first assertion of the lemma.
  By the first assertion of Lemma \ref{Greenberg compare lemma}, we have an exact sequence
\[0\lra N \lra X^{Gr}(A/F_{\infty}) \lra X(A/F_{\infty})\lra 0,\]
where $N$ is a finitely generated $\Op\ps{H}$-module. For notational simplicity, we write $X= X(A/F_{\infty})$ and $X^{Gr}=X^{Gr}(A/F_{\infty})$. For a sufficiently large $n$, we have $Z(\pi) = Z[\pi^n]$ for $Z = N, X^{Gr}, X$. We then have an exact sequence
\[ 0\lra N(\pi) \lra X^{Gr}(\pi) \lra X(\pi) \lra N/\pi^n. \]
Write $B=\im\big(X^{Gr}(\pi) \lra X(\pi)\big)$ and $C= \im\big(X(\pi) \lra N/\pi^n\big)$.  Consider the following two commutative diagrams
\[ \SelectTips{eu}{}
\xymatrix{
  0 \ar[r] & N(\pi) \ar[d] \ar[r]^{} & X^{Gr}(\pi)\ar[d]
  \ar[r]^{} & B\ar[d] \ar[r]& 0 \\
  0 \ar[r] & N \ar[r] & X^{Gr}
  \ar[r] & X \ar[r] &0   }
\]
\[ \SelectTips{eu}{}
\xymatrix{
  0 \ar[r] & B \ar[d] \ar[r]^{} & X(\pi) \ar[d]
  \ar[r]^{} & C \ar[r]& 0 \\
  & X  \ar@{=}[r]  & X
   &  &  }
\] with exact rows, where the vertical maps are the inclusion maps. From these two diagrams, we obtain two exact
sequences
\[ \ba{c}
 0\lra N_f \lra X^{Gr}_f \lra X/B \lra 0 \\
 0\lra C \lra X/B \lra X_f\lra 0.
 \ea \]
Since $N$ is finitely generated over $\Op\ps{H}$, it follows that $N_f$ and $C$ are finitely generated over $\Op\ps{H}$. Hence it follows from the two exact sequences that $X^{Gr}_f$ is finitely generated over $\Op\ps{H}$ if and only if $X_f$ is finitely generated over $\Op\ps{H}$. The second assertion of the lemma is then immediate from analysing the two exact sequences in the previous paragraph and noting the second assertion of Lemma
\ref{Greenberg compare lemma}.
\epf

\subsection{Selmer complexes}
\label{Selmer complexes}

 The notion of a Selmer
complex was first conceived and introduced in \cite{Nek}. In our
discussion, we consider a modified version of the Selmer complex as
given in \cite[4.2.11]{FK} (also see \cite[Section 3.1]{B}). There are several advantages
for using Selmer complexes in the formulation of the main conjecture rather than
the Pontryagin dual of the Selmer group. We will just mention two of these as word of mouth (for more details, see \cite{B, BV, FK, Nek}). Firstly, Selmer complexes have better functorial properties
than Selmer groups, and secondly, Nekov\'a\v{r} has shown that the Selmer complex is
able to explain certain phenomena concerning about the trivial zeroes (see \cite[0.10]{Nek}).

We now give the definition of the Selmer complex associated to the datum $\big(A,
\{A_v\}_{v|p}, \{A^+_v\}_{v|\R}\big)$. Write $T^* =
\Hom_{\cts}(A,\mu_{p^{\infty}})$ and $T^*_v =
\Hom_{\cts}(A/A_v,\mu_{p^{\infty}})$. For every finite extension $L$
of $F$ and $w$ a prime of $L$ above $p$, write $T^*_w = T^*_v$,
where $v$ is the prime of $F$ below $w$. For any profinite group
$\mathcal{G}$ and a topological abelian group $M$ with a continuous
$\mathcal{G}$-action, we denote by $C(\mathcal{G}, M)$ the complex
of continuous cochains of $\mathcal{G}$ with coefficients in $M$.
Let $F_{\infty}$ be an $S$-admissible extension of $F$ with Galois
group $G$. We define a $(\Op\ps{G})[G_S(F)]$-module
$\mathcal{F}_G(T^*)$ as follows: as an $\Op$-module,
$\mathcal{F}_G(T^*) = \Op\ps{G}\ot_{\Op}T^*$, and the action of
$G_S(F)$ is given by the formula $\sigma(x\ot t) =
x\bar{\sigma}^{-1}\ot \sigma t$, where $\bar{\sigma}$ is the
canonical image of $\sigma$ in $G \subseteq \Op\ps{G}$. We define
the $(\Op\ps{G})[\Gal(\bar{F}_v/F_v)]$-module $\mathcal{F}_G(T_v^*)$
in a similar fashion.

For every prime $v$ of $F$, we write $C\big(F_v,
\mathcal{F}_G(T^*)\big)= C\big(\Gal(\bar{F}_v/F_v),
\mathcal{F}_G(T^*)\big)$. For each prime $v$ not dividing $p$,
denote $C_f\big(F_v, \mathcal{F}_G(T^*)\big)$ to be the subcomplex
of $C\big(F_v, \mathcal{F}_G(T^*)\big)$, whose degree $m$-component
is $0$ unless $m\neq  0, 1$, whose degree $0$-component is
$C^0\big(F_v, \mathcal{F}_G(T^*)\big)$, and whose degree
$1$-component is
\[ \ker\Big( C^1\big(F_v, \mathcal{F}_G(T^*)\big)_{d=0} \longrightarrow
H^1\big(F_v^{ur},\mathcal{F}_G(T^*)\big) \Big). \]

The Selmer complex $SC(T^*, T^*_v)$ is then defined to be
\[ \mathrm{Cone}\bigg( C\big(G_S(F), \mathcal{F}_G(T^*)\big) \longrightarrow \hspace{3.5in}\]
\[ \hspace{0.5in} \longrightarrow
\displaystyle\bigoplus_{v|p}C\big(F_v,
\mathcal{F}_G(T^*)/\mathcal{F}_G(T^*_v)\big)\oplus \displaystyle
\bigoplus_{v\nmid p} C\big(F_v, \mathcal{F}_G(T^*)\big)/C_f\big(F_v,
\mathcal{F}_G(T^*)\big)\bigg)[-1].
\]
 Here $[-1]$ is the translation by $-1$ of the complex. We
write  $H^i\big(SC(T^*, T^*_v)\big)$ for the $i$th cohomology group
of the Selmer complex $SC(T^*, T^*_v)$. The Selmer complex, or
more precisely, its cohomology is related to other classical groups
arising in Iwasawa theory. For instance, if we have imposed trivial local condition
in the definition of the Selmer complex, then the cohomology of corresponding Selmer complex is
just the usual cohomology with compact support (see \cite[Subsection 1.6]{FK}, \cite[Subsection 2.3]{Kak}, \cite[Sections 4-5]{LimPT}, \cite[Section 4]{LimSh}, \cite[Chapter 5]{Nek}). In the general setting, the relation between the cohomology of the Selmer complexes and the various Iwasawa modules are extensively well-documented in
\cite[Chapters 6-9]{Nek} (also see \cite[Section 3]{B} or \cite[Subsection 4.2]{FK}). Interested readers are referred to these references for details. For the purpose of this paper, we just require
the following proposition which describes a fundamental relationship
between the Selmer complex and the Selmer groups. For its proof, we
shall refer readers to \cite[Proposition 4.2.35]{FK}.

\bp \label{FK prop} Let $\mathcal{G}$ be the kernel of
$\Gal(\bar{F}/F) \lra G$. For a place $v$ of $F$, fixing an
embedding $F \hookrightarrow F_v$, let $\mathcal{G}(v)$ be the
kernel of $\Gal(\bar{F}_v/F_v)\lra G$ and let $G_v\subseteq G$ be
the image. Then the following statements hold.
\begin{enumerate}
\item[$(a)$] $H^i\big(SC(T^*, T^*_v)\big) = 0$ for $i\neq 1, 2,
3$.

\item[$(b)$] We have an exact sequence
\[ \ba{c}
0 \lra X(A/F_{\infty}) \lra H^2\big(SC(T^*, T^*_v)\big) \lra
\displaystyle \bigoplus_{v|p} \Op\ps{G}\ot_{\Op\ps{G_v}}
\big(T^*_v(-1)\big)_{\mathcal{G}(v)}  \hspace{.5in}
 \\
\hspace{1.5in}\lra \big(T^*(-1)\big)_{\mathcal{G}} \lra H^3\big(SC(T^*,
 T^*_v)\big) \lra 0.
 \ea\]
\end{enumerate}
\ep

 In the next lemma, we write
$H^2\big(SC(T^*, T^*_v)\big)_f$ for the $\pi$-free quotient of
$H^2\big(SC(T^*, T^*_v)\big)$.

\bl \label{Greenberg equal complex}
  $X(A/F_{\infty})$ satisfies the $\M_H(G)$-property if and only if $H^2\big(SC(T^*, T^*_v)\big)$
satisfies the $\M_H(G)$-property.

Suppose further that $G=\Gal(F_{\infty}/F)$ is pro-$p$ without $p$-torsion and that for each $v|p$, the
decomposition group of $G$ at $v$ has dimension
$\geq 2$. Then if $X(A/F_{\infty})$ satisfies the $\M_H(G)$-property,
we have
\[ \rank_{\Op\ps{H}}(X_f(A/F_{\infty})) = \rank_{\Op\ps{H}}\Big(H^2\big(SC(T^*, T^*_v)\big)_f\Big).\]
 \el

\bpf
 Since $F_{\infty}$ contains $F_{\cyc}$, it follows that for every
prime $v|p$, the group $G_v$ has dimension at least 1.  Therefore,
$\bigoplus_{v|p} \Op\ps{G}\ot_{\Op\ps{G_v}}
\big(T^*_v(-1)\big)_{\mathcal{G}(v)}$ is finitely generated over
$\Op\ps{H}$. Now if for each $v|p$, the
decomposition group of $G$ at $v$ has dimension
$\geq 2$, then this is in fact a finitely generated torsion $\Op\ps{H}$-module. The conclusion of the lemma
now follows from a similar argument to that in Lemma \ref{Greenberg equal}.
\epf

\subsection{Some further remarks} \label{Some further remarks}

For the purposes of the formulation of the main conjecture, the Greenberg Selmer groups and the Selmer complexes are more natural algebraic objects to consider (see \cite{B,BV, FK,Gr89,Gr94,Kak,VBSD}). The aim of this paper is to study the invariance of the $\M_H(G)$-property of the Greenberg Selmer groups and the second cohomology groups of the Selmer complexes under congruence. In view of Lemmas
\ref{Greenberg equal} and \ref{Greenberg equal complex}, we are able to reduce the problem to considering the (strict) Selmer groups. When presenting and proving our main theorems in Section \ref{main results}, we shall always work with
$X(A/F_{\infty})$. Finally, the interested readers are invited to combine Lemmas \ref{Greenberg equal} and \ref{Greenberg equal complex} with \cite[Lemmas 3.1.2 and 3.2.2]{LimCMu} to see that $X(A/F_{\infty})$
captures much of the essential Iwasawa invariants (namely, the structure of $\pi$-primary submodules, the $\M_H(G)$-property and the $\Op\ps{H}$-rank of the $\pi$-free quotient) of the Greenberg Selmer groups and the Selmer complex $SC(T^*, T^*_v)$. We reassure the readers that this will be a profitable exercise.

\section{Main results} \label{main results}
In this section, we will prove the main theorem of the paper. We then apply our theorem to study the specialization of a big Galois representation. As mentioned in Subsection \ref{Some further remarks}, we only state formally the results for the strict Selmer groups, leaving the cases of the Greenberg Selmer groups and Selmer complexes for the readers to fill in.

\subsection{Congruent Galois representations} \label{Congruent Galois
representations}

As before, let $p$ be a prime. We let $F$ be a number
field. If $p=2$, we assume further that $F$ has no real primes.
Let $\big(A, \{A_v\}_{v|p}, \{A^+_v\}_{v|\R}
\big)$ and $\big(B, \{B_v\}_{v|p}, \{B^+_v\}_{v|\R}\big)$ be two data which
satisfy the conditions (a)--(d) as in Section 3. From now on, $S$ will always
denote a finite set of primes of $F$ which contains all the primes above
$p$, the ramified primes of $A$ and $B$, and the archimedean primes.
As before, for a given $S$-admissible $p$-adic Lie-extension $F_{\infty}$ of $F$, we
write $G = \Gal(F_{\infty}/F)$ and $H= \Gal(F_{\infty}/F^{\cyc})$.
\textit{We will always assume that the data obtained by base
change over every finite extension of $L$ of $F$ in $F_{\infty}$
also satisfy the conditions (a)--(d)}. Note that this assumption is automatically satisfied when $F_{\infty}$ is a pro-$p$ extension of $F$ by Lemma \ref{data base change}  (and noting our standing
assumption that if $p=2$, then $F$ has no real primes).

We now introduce the following important
congruence condition on $A$ and $B$ which allows us to be able to
compare the Selmer groups of $A$ and $B$.

\smallskip \noindent $\mathbf{(Cong_n)}$ : There is an isomorphism
$A[\pi^{n}]\cong B[\pi^{n}]$ of $G_S(F)$-modules which induces a
$\Gal(\bar{F}_v/F_v)$-isomorphism $A_v[\pi^{n}]\cong B_v[\pi^{n}]$
for every $v|p$.

For the purposes of comparing our results with that in \cite{LimCMu}, we first record the relevant
result (see \cite[Proposition 5.1.3]{LimCMu}).

\bp \label{MHG congruent} Let $F_{\infty}$ be an admissible $p$-adic
Lie extension of $F$, whose Galois group is a pro-$p$ group of
dimension 2 and has no elements of order $p$. Assume that
$A(F^{\cyc})$ and $B(F^{\cyc})$ are finite.  Suppose that
$\mathbf{(Cong_{\theta+1})}$ holds, where $\theta =
\max\{\theta_{\Op\ps{\Ga}}((X(A/F^{\cyc})), \theta_{\Op\ps{G}}((X(A/F_{\infty}))\}$. Then if
$X(A/F_{\infty})$ satisfies the $\M_H(G)$-property, so does
$X(B/F_{\infty})$. \ep

\smallskip
We may now present the main theorem of this paper. To simplify notation, we shall write $e_G(A) =
e_{\Op\ps{G}}\big(X(A/F_{\infty})\big)$. The following is the our theorem.

\bt \label{congruent thm} Let $F_{\infty}$ be an $S$-admissible $p$-adic Lie
extension of $F$. Suppose that the following statements hold.

\begin{enumerate}
\item[$(1)$] $\mathbf{(Cong_{e_G(A)+1})}$ holds.
\item[$(2)$] $X(A/F_{\infty})$ satisfies the $\M_H(G)$-property.
\item[$(3)$] $X(B/F_{\infty})$ has no nonzero pseudo-null $\Op\ps{G}$-submodules.
\end{enumerate}

Then we have the following statements.
\begin{enumerate}
\item[$(a)$] $X(B/F_{\infty})$ satisfies the $\M_H(G)$-property.
\item[$(b)$] Suppose further that the following statements hold.
  \begin{enumerate}
   \item[$(i)$] $F_{\infty}$ is a strongly $S$-admissible pro-$p$ $p$-adic Lie extension of $F$.
   \item[$(ii)$] For every $v\in S$, the decomposition group of $G$ at $v$ has dimension $\geq 2$.
   \end{enumerate}
    Then we have
   \[ \rank_{\Op\ps{H}}\big(X_f(A/F_{\infty})\big) = \rank_{\Op\ps{H}}\big(X_f(B/F_{\infty})\big).\]
   \end{enumerate} \et

\br
Readers will have observed that the assumptions of Theorem \ref{congruent thm} differs from Proposition \ref{MHG congruent}. This is due to the different approaches we have adopted in two proofs (which we have mentioned in the introductory section). The approach we used here enables us to obtain the preservation of the $\M_H(G)$-property for a general $p$-adic Lie extension, and at the same time, allows us to compare the $\Op\ps{H}$-rank of the $\pi$-free quotient of the dual Selmer groups. Of course, readers will have observed that this comes with some cost which we now explain. As mentioned in the introductory section, our approach follows that in \cite{AS}. There it is required that the dual Selmer group has no finite submodule, and therefore, in transferring their approach to the situation of a general admissible $p$-adic Lie extension, we are forced to work under assumption (3). The validity of Assumption (3) is known in quite a number of situations for Selmer groups attached to elliptic curves with good ordinary reduction at the prime $p$ (see \cite{HO, HV, Ochi, OcV02}). Thankfully, these observations can be partially extended to the general setting and this will be discussed in Appendix \ref{pseudo-null section}. The assumptions in (b)(i) and (ii) are rather mild. For instance, it is well-known that many admissible $p$-adic Lie extensions of interest satisfy assumption (b)(ii) (see \cite[Lemma 2.8(ii)]{C99} and \cite[Lemma 3.9]{HV}).

In Subsection \ref{big galois compare section}, we shall see how assumption (1) can be satisfied in big Galois representations. Of course, it will be of interest to come up with examples that do not come from big Galois representations. A natural source to start to look at will be \cite{AS, BS}. However, the examples in these cited works only consider the cyclotomic Iwasasa invariants over $\Q^{\cyc}/\Q$. But we do not know how the Iwasawa $\mu_{\Zp\ps{\Ga}}$-invariant varies upon base changing to $\Q(\mu_p)^{\cyc}/\Q(\mu_p)$. Secondly, although under the validity of $\M_H(G)$-property, the $\mu$-invariant has nice descent properties (see \cite[Lemmas 5.3 and 5.4]{CFKSV}, \cite[Corollary 3.2]{CS12} or \cite[Theorem 3.1]{LimMHG}), it is a mystery how the $\theta$-invariant (or the $e$-invariant) behaves under descent. At this point of writing, the author does not know how to resolve these two issues. Hence we shall content with ourselves with applying our theorem to specializations of a big Galois representation in this paper. We should also mention the second question also
implicitly comes up in \cite{LimCMu}, as the readers can see from the (rather imprecise) congruence condition imposed in Proposition \cite[Proposition 5.1.3]{LimCMu}.
\er

\medskip
In preparation of the
proof, we first introduce the ``mod $\pi^n$" Selmer group which is a standard object to work with in the study
of Selmer groups of congruent Galois representations (for instance, see \cite{AS,BS,Ch09,EPW,GV,Ha,LimAk,LimCMu,Sh,SS,We}). Let $n\geq 1$. For every
finite extension $\mathcal{F}$ of $F^{\cyc}$, we define
$J_v(A[\pi^n]/\mathcal{F})$ to be
 \[  \bigoplus_{w|v}H^1(\mathcal{F}_w, A[\pi^n])~\mbox{or}~
  \bigoplus_{w|v}H^1(\mathcal{F}_w, A/A_v[\pi^n])\]
according as $v$ does not or does divide $p$. We then define
\[ J_v(A[\pi^n]/F_{\infty}) = \ilim_{\mathcal{F}} J_v(A[\pi^n]/\mathcal{F}),\]
where the direct limit is taken over all finite extensions
$\mathcal{F}$ of $F^{\cyc}$ contained in $F_{\infty}$. The mod $\pi^n$
Selmer group is then defined by
 \[ S(A[\pi^n]/F_{\infty}) = \ker\Big(H^1(G_S(F_{\infty}), A[\pi^n])
 \lra \bigoplus_{v\in S} J_v(A[\pi^n]/F_{\infty})
 \Big).\]
We denote by $X(A[\pi^n]/F_{\infty})$ the Pontryagin dual of
$S(A[\pi^n]/F_{\infty})$.

We introduce some more notation which will be needed in our discussion.
 Write $C_v =A$ for $v\nmid p$ and $C_v =
A/A_v$ for $v|p$. For every prime $w$ of $F_{\infty}$ above $v$, we write
$C_v(F_{\infty,w})= (C_v)^{\Gal(\bar{F}_v/F_{\infty,w})}$.
We are now in the position to prove
Theorem \ref{congruent thm}.

\bpf[Proof of Theorem \ref{congruent thm}]
 Firstly, we note that if
$G_0$ is an open normal subgroup of $G$ and $H_0=H\cap G_0$, then a
finitely generated torsion $\Op\ps{G}$-module satisfies
$\M_H(G)$-property if and only if it satisfies
$\M_{H_0}(G_0)$-property. Now by the theorem of Lazard (cf.\
\cite[Cor.\ 8.34]{DSMS}), a compact $p$-adic Lie group always
contains a open normal uniform pro-$p$ subgroup. Hence to prove part
(a) of the theorem, we may, and we will, assume that $F_{\infty}$ is
a strongly admissible pro-$p$ $p$-adic Lie extension of $F$. In other words, $G$ (and
$H$) is pro-$p$ without $p$-torsion. Write $e = e_G(A)$ and consider
the following diagram
\[  \entrymodifiers={!! <0pt, .8ex>+} \SelectTips{eu}{}\xymatrix{
    0 \ar[r]^{} & S(A[\pi^{e+1}]/F_{\infty}) \ar[d] \ar[r] &
    H^1(G_S(F_{\infty}), A[\pi^{e+1}])
    \ar[d]_{b}
    \ar[r]^{} & \bigoplus_{v\in S}J_v(A[\pi^{e+1}]/F_{\infty}) \ar[d]_{c} \\
    0 \ar[r]^{} & S(A/F_{\infty})[\pi^{e+1}] \ar[r]^{}
    & H^1(G_S(F_{\infty}), A)[\pi^{e+1}] \ar[r] & \
    \bigoplus_{_{v\in S}}J_v(A/F_{\infty})[\pi^{e+1}]  } \]
with exact rows. We write $c =\oplus_w c_w$, where $w$ runs over the
set of primes of $F^{\cyc}$ above $S$. Denote $H_w$ to be the
decomposition group of $F_{\infty}/F^{\cyc}$ corresponding to a
fixed prime of $F_{\infty}$, which we also denote by $w$, above $w$.
Then we have $\ker c_w =
\mathrm{Coind}^{H_w}_H\big(C_v(F_{\infty,w})/\pi^{e+1}\big)$, where
$v$ is the prime of $F$ below $w$. This is clearly finitely
generated over $\Op\ps{H}$ and is annihilated by $\pi^{e+1}$. On the other hand, the long exact sequence in cohomology arising from
$0\lra A[\pi^{e+1}]\lra A\lra A\lra 0$ shows that the map $b$ is surjective with $\ker b
= A(F_{\infty})/\pi^{e+1}$ being finite. Hence we have an exact
sequence
\[ 0 \lra C\lra X(A/F_{\infty})/\pi^{e+1} \lra  X(A[\pi^{e+1}]/F_{\infty})\lra D \lra 0 \]
of $\Op\ps{H}$-modules which are annihilated by $\pi^{e+1}$.
Here $C$ is a subquotient of the Pontryagin dual of $\ker c$ and $D$
is a quotient of the Pontryagin dual of $\ker b$. In particular, $C$
is a finitely generated $\Op\ps{H}$-module and $D$ is finite. Therefore, we may apply Proposition \ref{alg compare} to conclude that $\frac{X(A[\pi^{e+1}]/F_{\infty})}{X(A[\pi^{e+1}]/F_{\infty})[\pi^{e}]}$
is finitely generated over $k\ps{H}$ if and only if
$\frac{X(A/F_{\infty})/\pi^{e+1} }{\left(X(A/F_{\infty})/\pi^{e+1}\right)[\pi^{e}]}$ is finitely generated over $k\ps{H}$. By
Proposition \ref{alg quotient compare}, the latter is isomorphic to
$X_f(A/F_{\infty})/\pi$ and this module is in turn finitely generated over $k\ps{H}$ as a consequence of the hypothesis that $X(A/F_{\infty})$
satisfies $\M_H(G)$-property.

Now it follows from $\mathbf{(Cong_{e_G(A)+1})}$ that there is an isomorphism
\[ S(A[\pi^{e+1}]/F_{\infty}) \cong S(B[\pi^{e+1}]/F_{\infty})\]
of $\Op\ps{G}$-modules (and hence $\Op\ps{H}$-modules) which in turn
induces an isomorphism
\[ \frac{X(A[\pi^{e+1}]/F_{\infty})}{X(A[\pi^{e+1}]/F_{\infty})[\pi^{e}]}
 \cong \frac{X(B[\pi^{e+1}]/F_{\infty})}{X(B[\pi^{e+1}]/F_{\infty})[\pi^{e}]}\]
of $k\ps{H}$-modules. Therefore, it follows that
$\frac{X(B[\pi^{e+1}]/F_{\infty})}{X(B[\pi^{e+1}]/F_{\infty})[\pi^{e}]}$
is finitely generated over $k\ps{H}$, and by Proposition \ref{alg quotient compare}, so is $X_f(B/F_{\infty})/\pi$. It then follows from an
application of Nakayama Lemma that $X_f(B/F_{\infty})$ is finitely
generated over $\Op\ps{H}$. This completes the proof of part (a)
of the theorem.

We now proceed with the proof of part (b)
of the theorem. By assumption (ii) and the fact that
$C_v(F_{\infty_w})/\pi^{e+1}$ is finite, we may apply Lemma
\ref{ind} to conclude that $(\ker c)^{\vee}/(\ker c)^{\vee}[\pi^e]$
and $(\ker c)^{\vee}/\pi^e$ are finitely generated torsion
$k\ps{H}$-modules. By Lemma \ref{useful lemma}, this in turn implies
that $C/C[\pi^e]$ and $C/\pi^e$ are finitely generated torsion
$k\ps{H}$-modules. Therefore, it follows from Propositions \ref{alg
compare} and \ref{alg quotient compare} that
\[ \ba{c} \rank_{k\ps{H}} \Big(X_f(A/F_{\infty})/\pi\Big) = \rank_{k\ps{H}} \Bigg(\displaystyle\frac{X(A[\pi^{e+1}]/F_{\infty})}{X(A[\pi^{e+1}]/F_{\infty})[\pi^{e}]}\Bigg) \\
\hspace{1.8in} = \rank_{k\ps{H}} \Bigg(\displaystyle\frac{X(B[\pi^{e+1}]/F_{\infty})}{X(B[\pi^{e+1}]/F_{\infty})[\pi^{e}]}\Bigg). \ea \]
Here the second equality follows from the above isomorphism deduced
from the $\mathbf{(Cong_{e_G(A)+1})}$ condition. By a similar argument
as above, we obtain an exact sequence
\[ 0 \lra C'\lra X(B/F_{\infty})/\pi^{e+1} \lra  X(B[\pi^{e+1}]/F_{\infty})\lra D' \lra 0 \]
of $\Op\ps{H}$-modules which are annihilated by $\pi^{e+1}$, where
$C'/C'[\pi^e]$ and $C'/\pi^e$ are finitely generated
torsion $k\ps{H}$-modules, and $D'$ is finite. By Proposition \ref{alg compare}, we then have
\[\ba{c} \rank_{k\ps{H}} \Big(X_f(A/F_{\infty})/\pi\Big) = \rank_{k\ps{H}} \Bigg(\displaystyle\frac{X(B[\pi^{e+1}]/F_{\infty})}{X(B[\pi^{e+1}]/F_{\infty})[\pi^{e}]}\Bigg) \\
 \hspace{1.8in} =
\rank_{k\ps{H}} \Bigg(\displaystyle\frac{X(B/F_{\infty})/\pi^{e+1}}{\big(X(B/F_{\infty})/\pi^{e+1}\big)[\pi^{e}]}\Bigg).  \ea\]
 But the last quantity is precisely $\rank_{k\ps{H}} \Big(X_f(B/F_{\infty})/\pi\Big)$
by Proposition \ref{upper bound}.
The required equality now follows an application of \cite[Corollary
1.10]{Ho} (or \cite[Proposition 4.12]{LimFine}) and the fact that
the $\pi$-free quotients of the dual Selmer groups have no
$\pi$-torsion. \epf

\subsection{Comparing specializations of a big Galois representation}
\label{big galois compare section}

We now apply the main result in Subsection \ref{Congruent Galois
representations} to compare the Selmer groups of specializations of
a big Galois representation. As before, let $p$ be a prime. We let
$F$ be a number field. If $p=2$, we assume further that $F$ has no
real primes. Denote $\Op$ to be the ring of integers of some finite
extension $K$ of $\Qp$. We write $R=\Op\ps{T}$ for the power series
ring in one variable. Suppose that we are given the following set of
data:

\begin{enumerate}
 \item[(a)] $\A$ is a
 cofree $R$-module of $R$-corank $d$ with a
continuous, $R$-linear $\Gal(\bar{F}/F)$-action which is unramified
outside a finite set of primes of $F$.

 \item[(b)] For each prime $v$ of $F$ above $p$, $\A_v$ is a
$\Gal(\bar{F}_v/F_v)$-submodule of $\A$ which is cofree of
$R$-corank $d_v$.

\item[(c)] For each real prime $v$ of $F$, we write $\A_v^+=
\A^{\Gal(\bar{F}_v/F_v)}$ which we assume to be cofree of $R$-corank
$d_v^+$.

\item[(d)] The quantities $d, d_v$ and $d_v^+$ satisfy the following identity
  \begin{equation} \label{data equality}
  \sum_{v|p} (d-d_v)[F_v:\Qp] = dr_2(F) +
 \sum_{v~\mathrm{real}}(d-d^+_v),
  \end{equation}
  where $r_2(F)$ is the number of complex primes of $F$.
 \end{enumerate}

For any prime element $f$ of $\Op\ps{T}$ such that $\Op\ps{T}/f$ is
a maximal order, then we can obtain a datum $\big(\A[f],
\{\A_v[f]\}_{v|p}, \{\A^+_v[f]\}_{v|\R}\big)$ in the sense of
Section \ref{Arithmetic Preliminaries}. The next lemma has a easy
proof which is left to reader.

\bl \label{big galois lemma}
 Let $f$ and $g$ be prime elements of $\Op\ps{T}$ with $\pi^n|(f-g)$
 such that $\Op\ps{T}/f$ and $\Op\ps{T}/g$ are maximal orders.
 Then $\A[f, \pi^n] = \A[g,\pi^n]$. One also has similar conclusions
 for $\A_v$ and $\A^+_v$.
\el

The next proposition compares the $\pi$-free quotient of the
dual Selmer groups of various specializations of a big Galois
representation. For a real number $x$, we denote $\lceil x \rceil$
to be the smallest integer not less than $x$.
We can now state and prove the following theorem which is
a refinement of \cite[Proposition 8.6]{SS}, \cite[Corollary
4.37]{B} and \cite[Proposition 5.2.3]{LimCMu}.

\bt
 \label{big galois compare}
 Let $F_{\infty}$ be an admissible $p$-adic Lie extension of $F$. Let $f$
be a prime element of $\Op\ps{T}$ such that $\Op':=\Op\ps{T}/f$ is a
maximal order. Set $A = \A[f]$ and suppose that $X(A/F_{\infty})$
satisfies the $\M_H(G)$-property. Set
\[ n : =
\Bigg\lceil\frac{e_{\Op'\ps{G}}\big(X(A/F_{\infty})\big)+1}{m}\Bigg\rceil,\]
where $m$ is the ramification index of $\Op'/\Op$. Suppose that $g$ is a
prime element of $\Op\ps{T}$ with $\pi^{n}|f-g$ such that
$\Op\ps{T}/g$ is isomorphic to $\Op'$ and $X(\A[g]/F_{\infty})$ has no nonzero pseudo-null $\Op\ps{G}$-submodule. Then we have that
$X(\A[g]/F_{\infty})$ satisfies the $\M_H(G)$-property.

Suppose further that the following statements hold.
  \begin{enumerate}
   \item[$(i)$] $F_{\infty}$ is a strongly $S$-admissible pro-$p$ $p$-adic Lie extension of $F$.
   \item[$(ii)$] For every $v\in S$, the decomposition group of $G$ at $v$ has dimension $\geq 2$.
   \end{enumerate}
       Then we also have
   \[ \rank_{\Op'\ps{H}}\big(X_f(A/F_{\infty})\big) = \rank_{\Op'\ps{H}}\big(X_f(\A[g]/F_{\infty})\big).\]
    \et

\bpf
  Let $g$ be a prime element of $\Op\ps{T}$ which satisfies the hypothesis
in the theorem. Let $\pi'$ be a prime element of $\Op'$ and write
$B= \A[g]$. It follows from Lemma \ref{big galois lemma} that there
is an isomorphism of $G_S(F)$-modules $A[\pi'^{mn}]\cong
B[\pi'^{mn}]$ which induces an isomorphism of
$\Gal(\bar{F}_v/F_v)$-modules $A_v[\pi'^{mn}]\cong B_v[\pi'^{mn}]$
for each prime $v$ of $F$ above $p$. By our hypothesis of $n$, we
have $mn \geq e_{\Op'\ps{G}}(A)+1$. In particular, the congruence
hypothesis $\mathbf{(C_{e_{\Op'\ps{G}}(A)+1})}$ holds for $A$ and
$\A[g]$. Hence the conclusion of the theorem is now immediate from
an application of Theorem \ref{congruent thm}. \epf

\br \label{SSremark}
 As noted in the Introduction, the above equality of $\Op\ps{H}$-ranks supports the prediction in
\cite[Theorem 8.8]{SS}. The proof in \cite[Theorem 8.8]{SS} rests on the (presumably) much stronger assumption that the Selmer group of the big Galois representation satisfies the $\M_H(G)$-property. However, at our current state of knowledge, it is not even clear if this latter assumption is a consequence of the assumption that the Selmer group of every specialization satisfies the $\M_H(G)$-property (but see \cite[Proposition 5.4]{CS12} and \cite[Proposition 8.6]{SS} for discussion in this direction). Finally, we mention that it would be interesting to prove such an equality for specializations which are congruent by a power smaller than the $\pi$-exponent of the dual Selmer group under the hypothesis of Theorem \ref{big galois compare}. At present, we do not know how to prove this equality, even with the added assumption that the Selmer group of every specialization satisfies the $\M_H(G)$-property. The point is that our methods here do not seem to be able to compare the Selmer groups when the congruence is smaller than the $\pi$-exponent. Note that the method of \cite[Theorem 8.8]{SS} breaks down here too, as we do not know whether the Selmer group of the big Galois representation satisfies the $\M_H(G)$-property, and so we cannot compare the $\Op\ps{H}$-rank of the $\pi$-free quotient of the dual Selmer groups of the specialization with the $R$-free quotient of the $R\ps{H}$-rank of the dual Selmer group of the big Galois representation.
\er

\section{Appendices}

\subsection{Nonexistence of pseudo-null submodules} \label{pseudo-null section}

As our main result (Theorem \ref{congruent thm}) requires that our Selmer groups to have no
nonzero pseudo-null submodules, we shall briefly discuss this property here. Namely, we give a sufficient criterion of this property.
As before, let $p$ be a prime and $F$ a number
field. If $p=2$, assume further that $F$ has no real primes.
Let $\big(A, \{A_v\}_{v|p}, \{A^+_v\}_{v|\R}
\big)$ be a datum which
satisfy the conditions (a)-(d) as in Section 3. Denote by $S$
a finite set of primes of $F$ which contains all the primes above
$p$, the ramified primes of $A$ and the archimedean primes.
\textit{We will also assume that the data obtained by base
change over every finite extension of $L$ of $F$ in $F_{\infty}$
also satisfy the conditions (a)--(d)}. Note that this assumption is automatically satisfied when $F_{\infty}$ is a pro-$p$ extension of $F$ by Lemma \ref{data base change} (and noting our standing
assumption that if $p=2$, then $F$ has no real primes).
We can now state the main result of this appendix.

\bp \label{pseudo-null submodule}
Let $F_{\infty}$ be a strongly admissible $p$-adic Lie
extension of $F$. Suppose that the following statements are valid.
\begin{enumerate}
 \item[$(i)$] $X(A/F_{\infty})$ is a torsion $\Op\ps{G}$-module.
 \item[$(ii)$] For every $v\in S$, the decomposition group of $G$ at $v$ has dimension $\geq 2$, and for those $v$ above $p$, the decomposition group of $G$ at $v$ has dimension $\geq 3$.
    \end{enumerate}
Then $X(A/F_{\infty})$ has no nonzero pseudo-null $\Op\ps{G}$-submodules.
 \ep

We should mention that the first proof of such results in the noncommutative situation goes back to the work
 of Ochi and Venjakob \cite{OcV02} (see also \cite{Ochi}). Later, in \cite{HO}, Hachimori and Ochiai gave a proof which simplified part of the original argument in \cite{OcV02}, and it is this latter approach we adopt here.

To prepare for the proof of the proposition, we recall the (Tate) dual data of $\big(A, \{A_v\}_{v|p}, \{A^+_v\}_{v|\R}
\big)$. For a $\Op$-module $N$,
we denote $T_{\pi}(N)$ to be its $\pi$-adic Tate module, i.e.,
$T_{\pi}(N) = \plim_i N[\pi^i]$.  We then set $A^* =
\Hom_{\cts}(T_{\pi}(A),\mu_{p^{\infty}})$, where $\mu_{p^{\infty}}$ denotes
the group of all $p$-power roots of unity. Similarly, for each $v|p$
(resp., $v$ real), we set $A^*_v=
\Hom_{\cts}(T_{\pi}(A/A_v),\mu_{p^{\infty}})$ (resp., $(A^*)^+_v=
\Hom_{\cts}(T_{\pi}(A/A^+_v),\mu_{p^{\infty}})$).  It is an easy
exercise to verify that $\big(A^*, \{A_v^*\}_{v|p},
\{(A^*)^+_v\}_{v|\R} \big)$ satisfies equality (\ref{data
equality}). The Selmer group attached to this dual data
is denoted by $S(A^*/F_{\infty})$, whose Pontryagin dual is denoted to be $X(A^*/F_{\infty})$.

\bl \label{surjective}
 Retain the assumptions in Proposition \ref{pseudo-null submodule}. Then we have
 $H^2(G_S(F_{\infty}),A)=0$ and a short exact sequence
 \[ 0 \lra S(A/F_{\infty})\lra H^1(G_S(F_{\infty}),A)\lra \bigoplus_{v\in S}J_v(A/F_{\infty})\lra 0.\]
 \el

\bpf
Since $X(A/F_{\infty})$ is a torsion $\Op\ps{G}$-module, it follows from \cite[Corollary 4.1.2]{LimCMu} that $X(A^*/F_{\infty})$ is also a torsion $\Op\ps{G}$-module. Recall that the fine Selmer group of $A^*$ (in the sense of \cite{CS}; see also \cite{LimFine}) is defined by the following exact sequence
\[
  0\lra R(A^*/F_{\infty})\lra H^1(G_S(F_{\infty}), A^*)
\lra
 \displaystyle\bigoplus_{v\in S}K_v(A^*/F_{\infty}),
  \]
  Here $K_v^i(A^*/F_{\infty}) = \ilim_L \bigoplus_{w|v}H^i(L_w, A^*)$, where the direct limit is taken over all finite extensions $L$ of $F$ contained in $F_{\infty}$ under the restriction maps. We write $Y(A^*/F_{\infty})$ for the Pontrayagin dual of $R(A^*/F_{\infty})$. Clearly, we have an injection $R(A^*/F_{\infty})\hookrightarrow S(A^*/F_{\infty})$ which in turn induces a surjection $X(A^*/F_{\infty})\twoheadrightarrow Y(A^*/F_{\infty})$. Since $X(A^*/F_{\infty})$ is torsion over $\Op\ps{G}$, so is $Y(A^*/F_{\infty})$. It then follows from \cite[Lemma 7.1]{LimFine} that $H^2(G_S(F_{\infty}),A)=0$ and this proves the first assertion of the lemma.

Now combining the Poitou-Tate sequence with the assertion of $H^2(G_S(F_{\infty}),A)=0$, we obtain
  an exact sequence
 \[
  0\lra S(A/F_{\infty})\lra H^1(G_S(F_{\infty}), A)
 \lra
 \displaystyle\bigoplus_{v\in S}J_v(A/F_{\infty})
  \lra
 \big(\widehat{S}(A^*/F_{\infty})\big)^{\vee}
 \lra 0
 \] and an injection
 \[ \plim_{L} H^1(G_S(L), T_{\pi} A^*) \lra \plim_{L}
 \bigoplus_{w\in S_L} H^1(L_w , T_{\pi}A^*), \]
  Here $\widehat{S}(A^*/F_{\infty})$ is defined as the kernel of the
map
 \[ \plim_{L} H^1(G_S(L), T_{\pi} A^*) \lra \plim_L
 \bigoplus_{w\in S_L}T_{\pi} H^1(L_w , A^*), \] where the inverse limit is
taken over all finite extensions $L$ of $F$ contained in $F_{\infty}$.
By \cite[Lemma 5.4]{OcV02}, $\plim_{L}
 \bigoplus_{w\in S_L} H^1(L_w , T_{\pi}A^*)$ is a torsionfree $\Op\ps{G}$-module. It then follow from this and the above injection that $\widehat{S}(A^*/F_{\infty})^{\vee}$ is also a torsionfree $\Op\ps{G}$-module.
On the other hand, taking the torsionness of $X(A/F_{\infty})$, the vanishing of $H^2(G_S(F_{\infty}, A)$, property (d) of our datum and the formulas in \cite[Theorem 4.1]{OcV03} and \cite[Theorem 7.1]{HV} into account, followed by a straightforward rank calculation, we have that $\widehat{S}(A^*/F_{\infty})^{\vee}$ has zero $\Op\ps{G}$-rank. Hence this forces $\widehat{S}(A^*/F_{\infty})=0$ and the required short exact sequence is a consequence of this.
\epf

We finally finish with the proof of the main result of this subsection.

\bpf[Proof of Proposition \ref{pseudo-null submodule}]
  By Lemma \ref{surjective}, we have a short exact sequence
  \[ 0 \lra \left(\bigoplus_{v\in S}J_v(A/F_{\infty})\right)^{\vee}\lra H^1(G_S(F_{\infty}),A)^{\vee}\lra X(A/F_{\infty})\lra 0.\]
The leftmost module in the sequence is a reflexive $\Op\ps{G}$-module by \cite[Lemma 5.4]{OcV02}. Since $H^2(G_S(F_{\infty}),A)=0$ by Lemma \ref{surjective}, it follows from an application of \cite[Theorem 4.7]{OcV02} that
$H^1(G_S(F_{\infty}),A)^{\vee}$ has no nonzero pseudo-null $\Op\ps{G}$-submodules. We may therefore now apply the criterion of Hachimori-Ochiai (cf. \cite[Proposition 3.5]{HO}) to conclude that $X(A/F_{\infty})$  has no nonzero pseudo-null $\Op\ps{G}$-submodules. \epf

\subsection{Comparing characteristic elements of Selmer groups}
\label{char subsection}

In this subsection, we aaply our main results to study the
characteristic elements of the Selmer groups. In preparation of this,
we need to recall
some further notion and notation from \cite{CFKSV}. Let
\[ \Si = \big\{\,s\in \Op\ps{G}~ \big| ~\Op\ps{G}/\Op\ps{G}s
~\mbox{is a finitely generated $\Op\ps{H}$-module} \big\}.\]
 By \cite[Theorem 2.4]{CFKSV}, $\Si$ is a left and right Ore set
consisting of non-zero divisors in $\Op\ps{G}$. Set $\Si^* =
\cup_{n\geq 0}\pi^n\Si$. It follows from \cite[Proposition
2.3]{CFKSV} that a finitely generated $\Op\ps{G}$-module $M$ is
annihilated by $\Si^*$ if and only if $M$ satisfies the $\M_H(G)$-property.
By abuse of notion, we shall denote $\M_H(G)$ to be the
category of all finitely generated $\Op\ps{G}$-modules which are
$\Si^*$-torsion.
It follows from the discussion in \cite[Section 3]{CFKSV} that we have the following
exact sequence
\[ K_1(\Op\ps{G}) \lra K_1(\Op\ps{G}_{\Si^*})
\stackrel{\partial_G}{\lra} K_0(\M_H(G))\lra 0 \] of $K$-groups. For
each $M$ in $\M_H(G)$, we define a \textit{characteristic element}
for $M$ to be any element $\xi_M$ in $K_1(\Zp\ps{G}_{\Si^*})$ which has
the property that
\[\partial_G(\xi_M) = -[M].\]

Let $\rho:G\lra GL_m(\Op_{\rho})$ denote a continuous group
representation, where
$\Op'= \Op_{\rho}$ is the ring of integers of some finite extension
of $K$. For $g\in G$, we write $\bar{g}$ for its image in $\Ga=G/H$.
We define a continuous group homomorphism
 \[G \lra M_d(\Op')\ot_{\Op}\Op\ps{\Ga}, \quad g\mapsto \rho(g)\ot
 \bar{g}. \] By \cite[Lemma 3.3]{CFKSV}, this in turn induces a map
 \[ \Phi_{\rho}: K_1(\Op\ps{G}_{\Si^*})\lra Q_{\Op'}(\Ga)^{\times}, \]
where $Q_{\Op'}(\Ga)$ is the field of fraction of $\Op'\ps{\Ga}$.
Let $\varphi: \Op'\ps{\Ga}\lra \Op'$ be the augmentation map and
denote its kernel by $\mathfrak{p}$. One can extend $\varphi$ to a
map $\varphi : \Op'\ps{\Ga}_{\mathfrak{p}}\lra K'$, where $K'$ is
the field of fraction of $\Op'$. Let $\xi$ be an arbitrary element
in $K_1(\Op\ps{G}_{\Si^*})$. If $\Phi_{\rho}(\xi)\in
\Op'\ps{\Ga}_{\mathfrak{p}}$, we define $\xi(\rho)$ to be
$\varphi(\Phi_{\rho}(\xi))$. If $\Phi_{\rho}(\xi)\notin
\Op'\ps{\Ga}_{\mathfrak{p}}$, we set $\xi(\rho)$ to be $\infty$.

Suppose for now  $G$ (and hence $H$) has no $p$-torsion. Following \cite{CSS}, we say that
the \textit{Akashi series} of $M$ exists if $H_i(H,M)$ is
$\Op\ps{\Ga}$-torsion for every $i$. In the case of this event, we
denote $Ak_H(M)$ to be the \textit{Akashi series} of $M$ which is
defined to be \[ \displaystyle \prod_{i\geq 0} g_i^{(-1)^i},\]
 where $g_i$ is the characteristic polynomial of $H_i(H,M)$.
 Of course, the Akashi series is only well-defined up to a
unit in $\Op\ps{\Ga}$. Also, note that since $G$ (and hence $H$) has
no $p$-torsion, $H$ has finite $p$-cohomological dimension (cf.\
\cite[Corollaire 1]{Ser}), and therefore, the alternating product is
a finite product. We can now state the following result which answer
\cite[Conjecture 4.8 Case 4]{CFKSV} partially, and is proven in
\cite[Proposition 6.2]{LimMHG} and \cite[Proposition 6.2]{LimAk}.

\bp \label{conjecture 4.8} Let $M\in \M_H(G)$. Suppose that $M$
contains no nonzero pseudo-null $\Op\ps{G}$-submodules. Let
$\xi_M$ be a characteristic element of $M$. Then the following
statements are equivalent.

\begin{enumerate}
\item[$(a)$] $\xi_M \in \al(\Op\ps{G}^{\times})$, where $\al$ is
the map $\Op\ps{G}_{\Si^*}^{\times}\lra K_1(\Op\ps{G}_{\Si^*})$.

\item[$(b)$] $\xi_M(\rho)$ is finite and lies in $\Op_{\rho}^{\times}$ for every continuous
group representation $\rho$ of $G$.

\item[$(c)$] $\Phi_{\rho}(\xi_M) \in \Op_{\rho}\ps{\Ga}^{\times}$ for every continuous
group representation $\rho$ of $G$.

\item[$(d)$] $\Phi_{\rho}(\xi_M) \in \Op_{\rho}\ps{\Ga}^{\times}$ for every Artin
 representation $\rho$ of $G$.

 \item[$(e)$] There exists an open normal pro-$p$ subgroup $G'$ of
$G$ such $Ak_{H'}\big(X(A/F_{\infty})\big) \in
\Op\ps{\Ga'}^{\times}$. Here $H' = H\cap G'$ and $\Ga' = G'/H'$.
 \end{enumerate}
\ep

\br
 By \cite[Lemma 4.9]{CFKSV}, the implications $(a)\Rightarrow(b)\Leftrightarrow(c)\Rightarrow(d)$ always hold without the pseudo-nullity assumption on $M$. By examining the proof of \cite[Proposition 6.2]{LimAk} (or rather \cite[Proposition 6.1]{LimAk}), we see that the implication $(d)\Rightarrow(e)$ also holds without the pseudo-nullity assumption. It is the implication $(e)\Rightarrow(a)$, where the pseudo-nullity assumption is required.
\er

We can now state and prove the following theorem which slightly refines
\cite[Theorem 6.3]{LimAk}. If $X(A/F_{\infty})$ satisfies the $\M_H(G)$-property, we denote by
$\xi_{A,f}$ a characteristic element of $X_f(A/F_{\infty})$. Now if $X(B/F_{\infty})$ satisfies the $\M_H(G)$-property, then
$\xi_{B,f}$ is defined similarly. We continue to write $\al_G$ for the natural
map
\[\Op\ps{G}_{\Si^*}^{\times}\lra K_1(\Op\ps{G}_{\Si^*}).\]

\bt \label{char congruent}
 Let $F_{\infty}$ be a strongly admissible $p$-adic Lie
extension of $F$. Suppose that the following conditions are satisfied.
\begin{enumerate}
 \item[$(a)$] The condition $\mathbf{(Cong_{\theta_G(A)+1})}$ is satisfied.
 \item[$(b)$] $X(A/F_{\infty})$ $($and hence $X(B/F_{\infty})$$)$ satisfies the $\M_H(G)$-property.
  \item[$(c)$] $X(B/F_{\infty})$ has no nonzero pseudo-null $\Op\ps{G}$-submodules.
  \item[$(d)$] For every $v\in S$, the decomposition group of $G$ at $v$ has dimension $\geq 2$.
  \end{enumerate}
Then if $\xi_{A,f}\in \al_G(\Op\ps{G}^{\times})$, so is $\xi_{B,f}$.
 \et

\bpf
  Suppose that $\xi_{A,f}\in \al_G(\Op\ps{G}^{\times})$. Then by
Proposition \ref{conjecture 4.8} and the remark thereafter, there exists an open normal
pro-$p$ subgroup $G'$ of $G$ such
$Ak_{H'}\big(X_f(A/F_{\infty})\big) \in \Op\ps{\Ga}^{\times}$, where
$H' = H\cap G'$. By \cite[Proposition 5.4]{LimMHG} or
\cite[Proposition 2.2]{LimAk}, this in turn implies that
$X_f(A/F_{\infty})$ is a finitely generated torsion
$\Op\ps{H'}$-module. By virtue of assumption (d) and Theorem
\ref{congruent thm}, it follows that $X_f(B/F_{\infty})$ is a
finitely generated torsion $\Op\ps{H'}$-module. By a well-known
theorem of Venjakob (cf. \cite[Example 2.3 and Proposition
5.4]{V03}), this in turn implies that $X_f(B/F_{\infty})$ is a
finitely generated pseudo-null $\Op\ps{G'}$-module. On the other
hand, it follows from assumption (c) and \cite[Lemma 4.2]{Su} that
$X_f(B/F_{\infty})$ has no nonzero pseudo-null
$\Op\ps{G}$-submodules. Consequently, we have $X_f(B/F_{\infty})=0$,
or $\partial_G(\xi_{B,f}) =0$. It then follows from the above exact
sequence of $K$-groups that there exists an element in
$K_1(\Op\ps{G})$ which maps to $\xi_{B,f}$. On the other hand, it is
well-known that $\Op\ps{G}^{\times}$ maps onto $K_1(\Op\ps{G})$, and
the required conclusion is now immediate from this. \epf

We end by recording the following result for completeness. Write $\xi_{A,\pi}$ (resp., $\xi_{B,\pi}$)
for a characteristic element of $X(A/F_{\infty})(\pi)$ (resp., a characteristic element of $X_f(B/F_{\infty})(\pi)$). Note that $\pi$-primary modules automatically satisfy $\M_H(G)$-property, and hence one can always attach characteristic elements to them.

\bt \label{char congruent}
 Let $F_{\infty}$ be a strongly admissible $p$-adic Lie
extension of $F$. Suppose that the following conditions are satisfied.
\begin{enumerate}
 \item[$(a)$] The condition $\mathbf{(Cong_{1})}$ is satisfied.
 \item[$(b)$] $X(A/F_{\infty})$ is a torsion $\Op\ps{G}$-module.
  \item[$(c)$]$\xi_{A,\pi}\in \al_G(\Op\ps{G}^{\times})$.
  \end{enumerate}
Then $X(B/F_{\infty})$ is a torsion $\Op\ps{G}$-module and $\xi_{B,\pi}\in \al_G(\Op\ps{G}^{\times})$.
 \et

\bpf
By Proposition \ref{conjecture 4.8} and the remark thereafter, there exists an open normal
pro-$p$ subgroup $G'$ of $G$ such
$Ak_{H'}\big(X(A/F_{\infty})(\pi)\big) \in \Op\ps{\Ga}^{\times}$, where
$H' = H\cap G'$. On the other hand, it is not difficult to deduce from the definition of the Akashi series and \cite[Lemma 2.2]{LimMHG} that
\[ Ak_{H'}\big(X(A/F_{\infty})(\pi)\big) = \pi^{\mu_{\Op\ps{G'}}\big(X(A/F_{\infty})\big)}.\]
Hence we conclude that
\[ \mu_{\Op\ps{G'}}\big(X(A/F_{\infty})\big) =0.\]
By \cite[Theorem 4.2.1]{LimCMu}, this in turn implies that $X(B/F_{\infty})$ is a torsion $\Op\ps{G}$-module
and \[ \mu_{\Op\ps{G'}}\big(X(B/F_{\infty})\big) =0.\]
By \cite[Definition 4.1]{BV}, $-1\in K_1(\Op\ps{G'}_{\Si})$ is a characteristic element for $X(B/F_{\infty})(\pi)$. Since $\Op\ps{G}^{\times}$ maps onto $K_1(\Op\ps{G})$, we have that
every characteristic element for $X(B/F_{\infty})(\pi)$ lies in $\al_{G'}(\Op\ps{G'}^{\times})$. It then follows from this and \cite[Theorem 6.8]{AW} that $\xi_{B,\pi}\in \al_G(\Op\ps{G}^{\times})$ which is what we want to show.
\epf

\begin{ack}
        The author would like to thank the anonymous referee for giving many valuable comments and suggestions
         that have improved the exposition of the paper. The author is supported by the
     National Natural Science Foundation of China under the Research
Fund for International Young Scientists
     (Grant No: 11550110172).
    \end{ack}

\end{document}